\documentstyle[11pt,epsfig,mst-stylefile,amssymb,amsmath,harvard]{article}
\newcommand{\bbZ}{{\Bbb Z}}
\newcommand{\bbR}{{\Bbb R}}
\newcommand{\bbN}{{\Bbb N}}
\newcommand{\bbC}{{\Bbb C}}

\newcommand{\norm}[1]{\left\|#1\right\|}

\renewcommand{\cite}{\citeyear}

\begin{document}

\title{Exponents, symmetry groups and classification of operator fractional Brownian motions
\thanks{The first author was supported in part by the Louisiana Board of Regents award LEQSF(2008-11)-RD-A-23.
The second author was supported in part by the NSF grants
DMS-0505628 and DMS-0608669.}
\thanks{The authors are thankful to Profs.\ Eric Renault and Murad Taqqu for their comments on this work.}
 }

\author{Gustavo Didier \footnote{Corresponding author.} \\ Tulane University \and Vladas Pipiras \\ CEMAT - Instituto Superior T\'{e}cnico and University of North Carolina}

\bibliographystyle{agsm}

\maketitle

\begin{abstract}
Operator fractional Brownian motions (OFBMs) are zero mean, operator self-similar (o.s.s.),
Gaussian processes with
stationary increments. They generalize univariate fractional
Brownian motions to the multivariate context. It is well-known that the so-called symmetry
group of an o.s.s.\ process is conjugate to subgroups of the orthogonal group.
Moreover, by a celebrated result of Hudson and Mason, the set of
all exponents of an operator self-similar process can be related to the tangent space
of its symmetry group.

In this paper, we revisit and study both the symmetry groups and
exponent sets for the class of OFBMs based on their spectral domain integral representations.
A general description of the symmetry groups of OFBMs
in terms of subsets of centralizers of the spectral domain parameters is provided.
OFBMs with symmetry groups of maximal
and minimal types are studied in any dimension. In particular,
it is shown that OFBMs have minimal symmetry groups (as thus,
unique exponents) in general, in the topological sense.
Finer classification results of OFBMs, based on the explicit construction of their symmetry groups,
are given in the lower dimensions $2$ and $3$. It is also shown
that the parametrization of spectral domain integral representations are, in a suitable sense,
not affected by the multiplicity of exponents, whereas the same is not true for
time domain integral representations.
\end{abstract}

\section{Introduction}
\label{s:intro}

This work is about the class of operator fractional Brownian
motions (OFBMs). Denoted by $B_H = \{B_H(t)\}_{t\in\bbR} =
\{(B_{H,1}(t),\ldots,B_{H,n}(t))' \in \bbR^n, t\in\bbR\}$, these
are multivariate zero mean Gaussian processes with stationary
increments which are operator self-similar (o.s.s.) with a matrix
exponent $H$. Operator self-similarity means that, for any
$c>0$,
\begin{equation}\label{e:oss}
    \{B_H(ct)\}_{t\in\bbR} \stackrel{{\mathcal L}}{=} \{c^H B_H(t)\}_{t \in \bbR},
\end{equation}
where $=_{\mathcal L}$ stands for the equality of finite-dimensional distributions and $c^H := e^{H\ln c}: =
\sum_{k=0}^\infty H^k (\ln c)^k/k!$. It is also assumed that OFBMs
are proper, in the sense that the support of the distribution of
$B_H(t)$ is $\bbR^n$ for every $t\in\bbR$. OFBMs play an important
role in the analysis of multivariate time series, analogous to
that of the usual fractional Brownian motion (FBM) in the univariate context. They have been
studied more systematically by Pitt \cite{pitt:1978}, Mason and Xiao
\cite{mason:xiao:2002}, Bahadaron, Benassi and D\c ebicki
\cite{bahadoran:benassi:debicki:2003}, Lavancier, Philippe and Surgailis
\cite{lavancier:philippe:surgailis:2009}, Didier and Pipiras
\cite{didier:pipiras:2009}, and others. Regarding
o.s.s.\ processes in general, see Hudson and Mason
\cite{hudson:mason:1982}, Laha and Rohatgi
\cite{laha:rohatgi:1981}, Sato \cite{sato:1991}, Maejima and Mason
\cite{maejima:mason:1994}, Maejima
\cite{maejima:1996,maejima:1998}, Meerschaert and
Scheffler \cite{meerschaert:scheffler:1999}, Section 11 in Meerschaert and Scheffler
\cite{meerschaert:scheffler:2001}, Chapter 9 in Embrechts and
Majima \cite{embrechts:maejima:2002}, Becker-Kern and Pap
\cite{becker-kern:pap:2008}. For related work on operator stable measures, see, for instance, Sharpe \cite{sharpe:1969}, Jurek and Mason
\cite{jurek:mason:1993}, Meerschaert and Veeh
\cite{meerschaert:veeh:1993,meerschaert:veeh:1995}, Hudson
and Mason \cite{hudson:mason:1981:JrMultAnalysis}, among others.

In particular, Didier and Pipiras \cite{didier:pipiras:2009} showed
that, under the mild assumption
\begin{equation}\label{e:eigen-assumption}
    0< \Re(h_k)<1,\quad k=1,\ldots,n,
\end{equation}
on the eigenvalues $h_k$ of the matrix exponent $H$, any OFBM
$B_H$ admits the so-called integral representation in the spectral
domain,
\begin{equation} \label{e:spectral-repres-OFBM}
\{B_{H}(t)\}_{t \in \bbR} \stackrel{{\mathcal L}}=
\Big\{\int_{\bbR} \frac{e^{itx} - 1}{ix}
(x^{-(H-\frac{1}{2}I)}_{+}A +
x^{-(H-\frac{1}{2}I)}_{-}\overline{A}) \widetilde{B}(dx)\Big\}_{t
\in \bbR}.
\end{equation}
Here, $x_{\pm} = \max\{\pm x,0\}$,
\begin{equation}\label{e:A}
    A = A_1 + i A_2
\end{equation}
is a complex-valued matrix with real-valued $A_1,A_2$,
$\overline{A}$ indicates the complex conjugate of $A$,
$\widetilde{B}(x) = \widetilde{B}_1(x) + i \widetilde{B}_{2}(x)$
is a complex-valued multivariate Brownian motion satisfying
$\widetilde{B}_1(-x) = \widetilde{B}_1(x)$, $\widetilde{B}_2(-x) =
-\widetilde{B}_2(x)$, and $\widetilde{B}_1$ and $\widetilde{B}_2$
are independent with induced random measure $\widetilde{B}(dx)$
satisfying $E\widetilde{B}(dx)\widetilde{B}(dx)^{*} = dx$. Thus,
according to (\ref{e:spectral-repres-OFBM}), OFBMs are
characterized (parametrized) by the matrices $H$ and $A$.

%In particular, there is no
%explicit identifiable parametrization available for OFBMs (i.e., a
%bijective map from the parameter space to the class of laws of
%OFBMs).

In this work, we continue the systematic study of OFBMs started in
Didier and Pipiras \cite{didier:pipiras:2009}. We now tackle the
issues of the symmetry structure of OFBMs and of the
non-uniqueness (multiplicity) of the exponents $H$, which, in our
view, are essentially unexplored. Such issues are strongly
connected. Since the fundamental work of Hudson and Mason
\cite{hudson:mason:1982}, it is well known that one given o.s.s.\
process $X$ may have multiple exponents. More specifically, if we
denote the set of exponents of $X$ by ${\cal E}(X)$, we have that
\begin{equation}\label{e:exponents}
{\mathcal E}(X) = H + T(G_X),
\end{equation}
where $H$ is any particular exponent of the process $X$. Here,
\begin{equation}\label{e:symmetry-group}
    G_X = \Big\{ C\in GL(n) : \{X(t)\}_{t\in \bbR} \stackrel{{\mathcal L}}{=}
    \{CX(t)\}_{t\in\bbR}     \Big\}
\end{equation}
is the so-called symmetry group of the process $X$ (where $GL(n)$
is the multiplicative group of invertible matrices), and
\begin{equation}\label{e:tangent-space}
    T(G_X) = \Big\{C: C = \lim_{n\to\infty} \frac{C_n - I}{d_n},\quad
    \mbox{for some}\ \{C_n\}\subseteq G_X,\ 0<d_n\to 0 \Big\}
\end{equation}
is the tangent space of the symmetry group $G_X$. By a result for
compact groups (e.g., Hoffman and Morris \cite{hoffman:morris:1998}, p.\ 49, or
Hudson and Mason \cite{hudson:mason:1982}, p.\ 285), it is known that
\begin{equation}\label{e:G_is_conjugate_to_compact_group}
G_{X} = W {\mathcal O}_{0}W^{-1}
\end{equation}
for some positive definite matrix $W$ and some subgroup ${\mathcal
O}_0$ of the orthogonal group. As a consequence, the knowledge
about (\ref{e:exponents}) is subordinated to that about the symmetry
group $G_X$ of $X$. For example, the exponent is unique for the
process $X$ if and only if the symmetry group $G_X$ is finite.

The description and study of symmetry groups beyond the
decomposition (\ref{e:G_is_conjugate_to_compact_group}) is a
reputedly difficult and interesting problem (see, for instance,
Billingsley \cite{billingsley:1966}, p.\ 176, and Jurek and Mason
\cite{jurek:mason:1993}, p.\ 60, both in the context of random
vectors; see also Meerschaert and Veeh \cite{meerschaert:veeh:1993,meerschaert:veeh:1995}). In this paper, we take up and provide some answers for
this challenging problem in the context of OFBMs. The main goal of this
paper is two-fold: to study the symmetry groups of OFBMs in as
much detail as possible, and based on this, to closely examine
(\ref{e:exponents}) for OFBMs $X = B_H$. We emphasize again
that, to the best of our knowledge, this is the first work where
symmetry groups are examined for any large class of o.s.s.\
processes (e.g., for the notion of symmetry groups of Markov processes, see
Liao \cite{liao:1992} and references therein). Indeed, since its
publication, the scope of the work of Hudson and Mason
\cite{hudson:mason:1982} appears to have remained only of general
nature, the same being true for the main result
(\ref{e:exponents}).

The integral representation (\ref{e:spectral-repres-OFBM})
provides a natural and probably the only means to consider
(almost) the whole class of OFBMs. Section \ref{s:groups} is
dedicated to the reinterpretation and explicit representation of
symmetry-related constructs in terms of the spectral
parametrization $H$, $A$. One of our main results provides a
decomposition of the symmetry groups of OFBMs into the
intersection of (subsets of) centralizers, i.e., sets of matrices
that commute with a given matrix. For example, in the case of time
reversible OFBMs (corresponding to the case when $AA^{*} =
\overline{AA^{*}}$), we show that the symmetry group $G_{B_H}$ is
conjugate to
\begin{equation}\label{e:symmetry-group-time-reversible}
    \bigcap_{x>0} G(\Pi_x).
\end{equation}
Here, $G(\Pi)$ denotes the centralizer of a matrix $\Pi$ in the
group $O(n)$ of orthogonal matrices, i.e.,
\begin{equation}\label{e:G-Pi}
    G(\Pi) = \{O\in O(n): O \Pi = \Pi O\},
\end{equation}
and the matrix-valued function $\Pi_x$ has the frequency $x$ as
the argument and is parametrized by $H$ and $A$. Moreover, which
is key for many technical results in this paper, we actually
express the positive definite conjugacy matrix $W$ in
(\ref{e:G_is_conjugate_to_compact_group}) in terms of the spectral
parametrization. This is a substantial improvement over previous
works on operator self-similarity, where only the existence of
such conjugacy is obtained, e.g., as in
(\ref{e:G_is_conjugate_to_compact_group}).

In view of (\ref{e:symmetry-group-time-reversible}) and
(\ref{e:G-Pi}), it is clear that the symmetry structure of OFBM is
rooted in centralizers. The characterization of the commutativity of
matrices is a well-studied algebraic problem (e.g., MacDuffee
\cite{macduffee:1946}, Taussky \cite{taussky:1953}, Gantmacher
\cite{gantmacher:1959}, Suprunenko and Tyshkevich
\cite{suprunenko:tyshkevich:1968}, Lax \cite{lax:2007}). We apply
the available techniques in a variety of ways to provide a
detailed study of the symmetry groups and the associated tangent
spaces (Sections \ref{s:min-max} and \ref{s:2-3}), as well as of
the consequences of the non-uniqueness (non-identifiability) of the parametrization for integral representations
(Section \ref{s:integral-repres}).

Our study of the symmetry structures of OFBMs is carried out from two
perspectives: first, by looking at the extremal cases, i.e.,
maximal and minimal symmetry for arbitrary dimension, and second, by
conveying a full description of \textit{all} symmetry groups in
the lower dimensions $n=2$ and $n=3$.

Section \ref{s:min-max} is dedicated to the first perspective. We
completely characterize OFBMs with maximal symmetry, i.e., those
whose symmetry groups are conjugate to $O(n)$. We establish the
general form of their covariance function and of their spectral
parametrization. However, as intuitively clear, maximal symmetry
corresponds to a strict subset of the parameter space of OFBMs. In
view of this, one can naturally ask what the most typical symmetry
structure for OFBMs is, in a suitable sense. A related question is
whether the multiplicity of exponents (and, thus, the
non-identifiability of the parametrization) is a general
phenomenon. Section \ref{s:min-max} contains our answer to both
questions, which is, indeed, one of our main results. We prove
that, in the topological sense, OFBMs with minimal symmetry groups
(i.e., $\{I,-I\}$) form the \textit{largest} class within all
OFBMs. As a consequence, in the same sense, OFBMs generally have
\textit{unique} (identifiable) exponents. To establish this result, in our
analysis of the centralizers $G(\Pi_x)$, we bypass the need to
deal with the major complexity of the eigenspace structure of the
function $\Pi_x$ by looking at its behavior at the origin of the
Lie group (i.e., as $x \rightarrow 1$), where a great deal of
information about $\Pi_x$ is available through the celebrated
Baker-Campbell-Hausdorff formula.

Section \ref{s:2-3} contains a full description of the symmetry
structure of low-dimensional OFBMs, namely, for dimensions $n=2$
and $n=3$. We provide a classification of OFBMs based on their
symmetry groups. For example, when $n=2$, the symmetry group of a
general OFBM can be, up to a conjugacy, of only one of the
following types:
\begin{itemize}
    \item [$(i)$] minimal: $\{I,-I\}$;
    \item [$(ii)$] trivial: $\{I,-I,R,-R\}$, where $R$ is a reflection matrix;
    \item [$(iii)$] rotational: $SO(2)$ (the group of rotation matrices);
    \item [$(iv)$] maximal: $O(2)$.
\end{itemize}
Such classification of types for $n=2$ stands in contrast with the
situation with random vectors, for which $SO(n)$ cannot be a
symmetry group (Billingsley \cite{billingsley:1966}).
Nevertheless, we show that the latter statement is \textit{almost}
true for OFBMs, since $SO(n)$ cannot be a symmetry group if $n
\geq 3$. In both $n=2$ and $n=3$, we provide examples of OFBMs in
all identified classes, and also discuss the structure of the
resulting exponent sets ${\cal E}(B_H)$.

In Section \ref{s:integral-repres}, we examine the consequences of
non-identifiability for integral representations of OFBMs. We show
that the multiplicity of the exponents $H$ does \textit{not}
affect the parameter $A$ in (\ref{e:spectral-repres-OFBM}) in the
sense that the latter can be chosen the same for any of the
exponents. Intriguingly, this turns out \textit{not} to be the
case for the parameters in the time domain representation of
OFBMs, and points to one advantage of spectral domain
representations.

%In face of the non-uniqueness of the parametrization of OFBMs, a
%natural pursuit is that of invariants and of an identifiable
%parametrization (see, for instance, Bahadaron et al.\
%\cite{bahadoran:benassi:debicki:2003}). This issue is explored separately in
%Didier and Pipiras \cite{didier:pipiras:2010:identif} in the context of general o.s.s.\ processes.

In summary, the structure of the paper is as follows. Some
preliminary remarks and notation can be found in Section
\ref{s:prelim}. Section \ref{s:groups} concerns structural results on
the symmetry groups of OFBMs. OFBMs with maximal and minimal
symmetry groups are studied in Section \ref{s:min-max}.
The classification of OFBMs according to their symmetry groups in the
lower dimensions $n=2$ and $n=3$ can be found in Section
\ref{s:2-3}. Section \ref{s:integral-repres} contains results on
the consequences of the non-uniqueness of the parametrization for integral
representations. The appendix contains several auxiliary facts
for the reader's convenience.

\section{Preliminaries}
\label{s:prelim}

\subsection{Notation}
We shall use throughout the paper the following notation for
finite-dimensional operators (matrices). All with respect to the
field $\bbR$, $M(n)$ or $M(n,\bbR)$ is the vector space of all $n
\times n$ operators (endomorphisms), $GL(n)$ or $GL(n,\bbR)$ is the
general linear group (invertible operators, or automorphisms),
$O(n)$ is the orthogonal group of operators $O$ such that $OO^{*} =
I = O^{*}O$ (i.e., the adjoint operator is the inverse), $SO(n)
\subseteq O(n)$ is the special orthogonal group of operators
(rotations) with determinant equal to 1, and $so(n)$ is the vector
space of skew-symmetric operators (i.e., $A^{*} = -A$). Similarly, $M(m,n,\bbR)$
is the space of $m \times n$ real matrices. The notation
will indicate the change to the field $\bbC$. For instance,
$M(n,\bbC)$ is the vector space of complex endomorphisms. Whenever
it is said that $A \in M(n)$ has a complex eigenvalue or eigenspace,
one is considering the operator embedding $M(n) \hookrightarrow
M(n,\bbC)$. $U(n)$ is the group of unitary matrices, i.e., $UU^*=I=U^*U$. ${\mathcal S}(n,\bbR)$ is the space of symmetric matrices.
 We will say that two endomorphisms $A, B \in M(n)$ are
\textit{conjugate} (or similar) when there exists $P \in GL(n,\bbC)$
such that $A = P B P^{-1}$. In this case, $P$ is called a
\textit{conjugacy}. The expression
$\textnormal{diag}(\lambda_1,\ldots,\lambda_n)$ denotes the
operator whose matrix expression has the values $\lambda_1,
\ldots, \lambda_n$ on the diagonal and zeros elsewhere. We make no
conceptual distinction between characteristic roots and
eigenvalues. We also write $S^{n-1}:=\{v \in \bbR^n:
\norm{v}=1\}$; in particular, we denote the complex sphere by $S^{2n-1}$. ${\mathbf 0}$ represents a matrix of zeroes of suitable dimension. Whenever necessary, we
will specify the dimension of the identity matrix by writing $I_n$. Unless otherwise stated, we consider the
so-called spectral matrix norm $\norm{\cdot}$, i.e., $\norm{A}$ is the square
root of the largest eigenvalue of $A^*A$. For $\{A_n\}_{n \in \bbN}$,
$A \in M(n,\bbC)$, we write $A_n \rightarrow A$ when $\norm{A_n
- A} \rightarrow 0$. $V^{\perp}$ is the subspace perpendicular to a given vector subspace $V$. For a set of (column) vectors $v_1,\hdots,v_n$, $A := (v_1,\hdots,v_n)$ is the matrix whose columns are such vectors. We denote the $i$-th Euclidean vector by $e_i$, $i=1,\hdots,n$.

Throughout the paper, we set
\begin{equation}\label{e:D}
    D = H - \frac{1}{2}I,
\end{equation}
for an operator exponent $H$. We shall also work with the real part
$\Re(AA^*) = A_1 A_1^* + A_2 A_2^*$ and the imaginary part
$\Im(AA^*)= A_2 A_1^* - A_1 A_2^*$ of $AA^*$. For the real part, in
particular, we will use the decomposition
\begin{equation}\label{e:decomp-Re-AA*}
    \Re(AA^*) = S_R \Lambda_R^2 S_R^* = W^2,
\end{equation}
with an orthogonal $S_R$, a diagonal $\Lambda_R$ and a positive
(semi-)definite
\begin{equation}\label{e:W}
    W = S_R \Lambda_R S_R^*.
\end{equation}
We shall use the assumption that
\begin{equation}\label{e:full-rank}
    \Re(AA^*)\ \mbox{has full rank},
\end{equation}
in which case $\Lambda_R$ in (\ref{e:decomp-Re-AA*}) has the
inverse $\Lambda_R^{-1}$. As shown in Didier and Pipiras
\cite{didier:pipiras:2009}, the condition (\ref{e:full-rank}) is
sufficient (though not necessary) for the integral in
(\ref{e:spectral-repres-OFBM}) to be proper and hence to define an
OFBM.

All through the paper, we assume $n \geq 2$.

\subsection{Remarks on the multiplicity of matrix exponents}

In this section, we make a few remarks to a reader less familiar
with the subject of this work. It may appear a bit surprising that
an o.s.s.\ process may have multiple exponents, as formalized in
(\ref{e:exponents}). This can be understood from at least two
inter-related perspectives: the properties of operator (matrix)
exponents and the distributional properties of o.s.s.\ processes.
From the first perspective, consider for example matrices of the
form
\begin{equation}\label{e:L_s}
L_s = \left(
\begin{array}{cc}
0 & s \\
-s & 0 \\
\end{array} \right) \in so(2),
\end{equation}
where $s \in \bbR$. Being normal, these matrices can be
diagonalized as $L_s = U_2 \Lambda_s U^{*}_2$, where $U_2 \in
U(2)$ and $\Lambda_s = \textnormal{diag}(is,-is)$. In particular,
$\exp(L_{2\pi k}) = I$, $k \in \bbZ$, since $e^{i2\pi k}=1$.
Since $L_s$ and $L_{s'}$ commute for any $s,s'\in \bbR$, this
yields
\begin{equation} \label{e:L_2.pi.k}
\exp{(L_s)} = \exp{(L_{2\pi k})} \exp{(L_s)} = \exp{(L_{2\pi k} +
L_s)},
\end{equation}
and shows the potential non-uniqueness of operator exponents
stemming from purely operator (matrix) properties. Note also that
the situation here is quite different from the 1-dimensional case:
in the latter, the same is possible but only with complex
exponents, whereas here the matrices $L_{2\pi k}$ have purely real
entries.

From the perspective of distributional properties, we can
illustrate several ideas through the following simple example. The
OFBMs in this example will appear again in Section \ref{s:min-max}
below.

\begin{example} (Single parameter OFBM)\label{ex:single-param_OFBM}
Consider an OFBM $B_H$ with covariance function $EB_H(t)B_H(s)^*
=: \Gamma(t,s) = \Gamma_{h}(t,s)I$, where $\Gamma_{h}(t,s)$ is the
covariance function of a standard univariate FBM with parameter $h
\in (0,1)$. This process is o.s.s.\ with exponent $H = h \hspace{0.5mm}I$, and
will be called a single parameter OFBM. Since $B_H$
is Gaussian, $O \in G_{B_H}$ if and only if $O \Gamma(t,s)O^{*} =
\Gamma(t,s)$. In the case of a single parameter OFBM, this is
equivalent to $OO^{*} = I$ or, since $O$ has an inverse ($B_H$ is
assumed proper), $OO^{*} = O^{*}O = I$. In other words, $G_{B_H} =
O(n)$ and
$$
{\mathcal E}(B_H) = H + T(O(n)) = H + so(n).
$$
Thus, a single parameter OFBM has multiple exponents. From
another angle, for a given $c > 0$ and $L \in so(n)$, we have
$L\log(c) \in so(n)$ and hence $\exp(L \log(c)) = c^{L} \in O(n) =
G_{B_H}$. Then,
$$
\{B_H(ct)\}_{t \in \bbR} \stackrel{{\mathcal L}}= \{c^{H}
B_H(t)\}_{t \in \bbR} \stackrel{{\mathcal L}}= \{c^{H}c^{L}
B_H(t)\}_{t \in \bbR} \stackrel{{\mathcal L}}= \{c^{H+L} B_H(t)\}_{t
\in \bbR},
$$
which also shows that the exponents are not unique.

For later use, we also note that an equivalent way to define a
single parameter OFBM is to say that it has the spectral
representation
\begin{equation}\label{e:spectral-repres-single-OFBM}
\{B_{H}(t)\}_{t \in \bbR} \stackrel{{\mathcal L}}= \Big\{C
\int_{\bbR} \frac{e^{itx} - 1}{ix} |x|^{-(h-\frac{1}{2})}
\widetilde{B}(dx)\Big\}_{t \in \bbR},
\end{equation}
where $C$ is an appropriate normalizing constant and $\widetilde
B(dx)$ is as in (\ref{e:spectral-repres-OFBM}).
\end{example}

\subsection{Basics of matrix commutativity}

We now recap some key facts and results about matrix commutativity
that are repeatedly used in the paper. To put it shortly, two matrices $A, B \in M(n,\bbC)$ commute if
and only if they share a common basis of generalized eigenvectors
(see Lax \cite{lax:2007}, p.\ 74). This means that there exists a
matrix $P \in GL(n,\bbC)$ such that we can write $A = P
J_{A}P^{-1}$ and $B = P J_{B}P^{-1}$, where $J_A$ and $J_B$ are in
Jordan canonical form. In particular, if $A$, $B$ are
diagonalizable, then they must share a basis of eigenvectors.
When, for example, $A = I$, we can interpret that $A$ commutes
with any $B = P J_{B} P^{-1} \in M(n,\bbC)$ because for (any) $P \in GL(n,\bbC)$, $A
= P I P^{-1}$.

A related issue is that of the characterization of the set of all
matrices that commute with a given matrix $A$, the so-called
centralizer ${\mathcal C}(A)$. In particular, one is often interested in
constructing the latter based on the Jordan decomposition of $A$.

Before enunciating the main theorem on ${\mathcal C}(A)$,
we look at an example adapted from Gantmacher \cite{gantmacher:1959}.

\begin{example} Assume the matrix $A \in M(10,\bbC)$, with Jordan representation $A = P J_{A} P^{-1}$,
has the elementary divisors (i.e., the characteristic polynomials of the Jordan blocks)
\begin{equation} \label{e:elem_divisors_A}
(\lambda - \lambda_1)^{3}, (\lambda - \lambda_1)^{2}, (\lambda -
\lambda_2)^{2}, (\lambda - \lambda_3),(\lambda - \lambda_3),(\lambda - \lambda_3),
\end{equation}
where the eigenvalues $\lambda_1, \lambda_2, \lambda_3$ are pairwise distinct. Then,
${\mathcal C}(A)$ consists of matrices of the form $X = P \widetilde{X} P^{-1}$, where
%\[ \left( \begin{array}{ccccccccccccc}
%a  &  0  & 0  & | & 0 & 0 & | & 0 & 0 & | & 0 & 0 & 0\\
%b  &  a  & 0  & | & f & 0 & | & 0 & 0 & | & 0 & 0 & 0\\
%c  &  b  & a  & | & g & f & | & 0 & 0 & | & 0 & 0 & 0\\
%-  &  -  & -  & | & - & - & | & - & - & | & - & - & -\\
%d  &  0  & 0  & | & h & 0 & | & 0 & 0 & | & 0 & 0 & 0\\
%e  &  d  & 0  & | & i & h & | & 0 & 0 & | & 0 & 0 & 0\\
%-  &  -  & -  & | & - & - & | & - & - & | & - & - & -\\
%0  &  0  & 0  & | & 0 & 0 & | & j & 0 & | & 0 & 0 & 0\\
%0  &  0  & 0  & | & 0 & 0 & | & k & j & | & 0 & 0 & 0 \\
%-  &  -  & -  & | & - & - & | & - & - & | & - & - & -\\
%0  &  0  & 0  & | & 0 & 0 & | & 0 & 0 & | & l & m & n \\
%0  &  0  & 0  & | & 0 & 0 & | & 0 & 0 & | & o & p & q \\
%0  &  0  & 0  & | & 0 & 0 & | & 0 & 0 & | & r & s & t \\
%\end{array} \right). \]\\
\begin{equation}\label{e:Xtilde}M(10,\bbC) \ni \widetilde{X} =  \left( \begin{array}{ccccccccccccccc}
a  &     &    & | &   &   & | &   &   & | &   & | &  & | &  \\
b  &  a  &    & | & f &   & | &   &   & | &   & | &  & | &  \\
c  &  b  & a  & | & g & f & | &   &   & | &   & | &  & | & \\
-  &  -  & -  & | & - & - & | & - & - & | & - & | &- & | &-\\
d  &     &    & | & h &   & | &   &   & | &   & | &  & | &\\
e  &  d  &    & | & i & h & | &   &   & | &   & | &  & | & \\
-  &  -  & -  & | & - & - & | & - & - & | & - & | &- & | &-\\
   &     &    & | &   &   & | & j &   & | &   & | &  & | & \\
   &     &    & | &   &   & | & k & j & | &   & | &  & | &  \\
-  &  -  & -  & | & - & - & | & - & - & | & - & | &- & | &-\\
   &     &    & | &   &   & | &   &   & | & l & | & m & | & n \\
-  &  -  & -  & | & - & - & | & - & - & | & - & | & - & | & -\\
   &     &    & | &   &   & | &   &   & | & o & | & p & | & q \\
-  &  -  & -  & | & - & - & | & - & - & | & - & | & - & | & -\\
   &     &    & | &   &   & | &   &   & | & r & | & s & | & t \\
\end{array} \right). \end{equation}
The blocks on the diagonal correspond to the Jordan blocks of $J_A$ and the empty entries are zeroes.
\end{example}

We now turn to the structure of the blocks for the general case, and (\ref{e:Xtilde}) serves as an illustration of the latter. We say that a matrix $X_{\alpha \beta} \in M(p_{\alpha},p_{\beta}, \bbC)$ has regular lower triangular form (e.g., as each of the blocks in (\ref{e:Xtilde}) with letters $a$ through $j$) if it can be written as
\begin{equation}\label{e:Xab}
X_{\alpha \beta}= \left\{\begin{array}{cc}
(T_{p_{\alpha}},{\mathbf 0}), & \textnormal{if }p_{\alpha} \leq p_{\beta},\\
({\mathbf 0}^{'} ,T^{'}_{p_{\beta}})^{'}, & \textnormal{if } p_{\alpha} > p_{\beta},\\
\end{array}\right.
\end{equation}
where $T_{p_{\alpha}} \in M(p_{\alpha},\bbC)$ is a Toeplitz lower triangular matrix. Also denote by $N_{p_{\alpha}} \in M(p_{\alpha}, \bbC)$
the nilpotent matrix
$$
N_{p_{\alpha}} =
\left(\begin{array}{cccc}
0 &        &          &      \\
1 &      0 &           &     \\
  & \ddots &  \ddots   &     \\
  &        &  1        &  0
\end{array}\right).
$$

The next theorem characterizes ${\mathcal C}(A)$. The proof can be found in Gantmacher \cite{gantmacher:1959}, p. 219 (see also pp.\ 220-224).

\begin{theorem} \label{t:comm_general}
Let $A \in M(n,\bbC)$, where $A = P J_{A} P^{-1}$ and $J_{A}$ is in Jordan canonical form, i.e.,
$$
J_{A} = \textnormal{diag}(\lambda_1 I_{p_1} + N_{p_1}, \hdots, \lambda_u I_{p_u} + N_{p_u})
$$
with (not necessarily distinct) eigenvalues $\lambda_1,\hdots,\lambda_u$. Then, the general solution to the equation $AX = XA$
is given by the formula $X = P X_{J_A} P^{-1}$, where $X_{J_A}$ is the general solution to the equation $J_A X_{J_A} = X_{J_A} J_A$. Here, $X_{J_A}$
can be decomposed into blocks $X_{\alpha \beta} \in M(p_{\alpha},p_{\beta},\bbC)$, $\alpha, \beta=1, \hdots,u$, where
$$X_{\alpha \beta} = \left\{\begin{array}{cccc}
{\mathbf 0}, & \textnormal{\textit{if }} \lambda_{\alpha} \neq \lambda_{\beta},\\
\textnormal{ \textit{as in (\ref{e:Xab})}}, & \textnormal{\textit{if }}\lambda_{\alpha} = \lambda_{\beta}.
\end{array}\right.
$$
\end{theorem}

\begin{corollary}\label{c:commute=>invariant}
Let $A \in {\mathcal S}(n,\bbR)$, where $A = O \Upsilon O^*$, $O \in O(n)$, $\Upsilon = \textnormal{diag}(\eta_1,\hdots,\eta_n)$. Assume $\eta_{j_1} = \hdots = \eta_{j_k}$ (possibly $k=1$) for some subset of eigenvalues of $A$, and $\eta_{j_1} \neq \eta_{i}$ for any other eigenvalue $\eta_i$ of $A$. If another matrix $M \in M(n,\bbR)$ commutes with $A$, then $M$ can be represented as
$$
M = (o_{j_1},\hdots,o_{j_k},o_{i_1},\hdots,o_{i_{n-k}})\textnormal{diag}(K_{11},K_{22})(o_{j_1},\hdots,o_{j_k},o_{i_1},\hdots,o_{i_{n-k}})^*,
$$
where $\textnormal{diag}(K_{11},K_{22})$ is block-diagonal with $K_{11} \in M(k,\bbR)$, $K_{22} \in M(n-k,\bbR)$, where $o_{j_1},\hdots,o_{j_k}$ are the column vectors in $O$ associated with the eigenvalues $\eta_{j_1}, \hdots, \eta_{j_k}$, respectively, and $o_{i_1},\hdots,o_{i_{n-k}}$ are the column vectors in $O$ associated with the eigenvalues $\eta_{i_1}, \hdots, \eta_{i_{n-k}}$, respectively.

Consequently, $\textnormal{span}_{\bbR}\{o_{j_1},\hdots,o_{j_k}\}$ is an invariant subspace with respect to $M$.
\end{corollary}

In view of Theorem \ref{t:comm_general} (and Corollary \ref{c:commute=>invariant}), it is intuitively clear
that, if a matrix $\Gamma$ commutes with two matrices $A$ and $B$
which exhibit completely different invariant subspaces, then
$\Gamma$ can only be a multiple of the identity. This is
accurately stated for the case of symmetric matrices in the next lemma, which is used in Section \ref{s:min-max}.

\begin{lemma}\label{l:irreducible}
Let
\begin{eqnarray}\label{e:Sneq}
A, B \in {\mathcal S}_{\neq} &:=& \{S \in {\mathcal S}(n,\bbR):\textnormal{$S$ has pairwise distinct eigenvalues}\}.
\end{eqnarray}
Assume $A$ and $B$ have no $k$-dimensional invariant subspaces in common for $k=1,\hdots,n-1$. If $M \in M(n,\bbR)$ commutes with both $A$ and $B$, then $M$ is a scalar matrix.
\end{lemma}

The proof of Corollary \ref{c:commute=>invariant} and Lemma \ref{l:irreducible} can be
found in Appendix \ref{s:auxiliary}, together with some additional
results on matrix commutativity and matrix representations.

\section{Symmetry groups of OFBMs}
\label{s:groups}

Consider an OFBM $B_H$ with the spectral representation
(\ref{e:spectral-repres-OFBM}). In this section, we provide some structural
results on the nature of the symmetry group $G_{B_H}$ (see (\ref{e:symmetry-group})). In particular, we
explicitly express it as an intersection of subsets of centralizers.

For notational simplicity, denote
$G_{B_H}$ by $G_H$. Since OFBMs are Gaussian and two Gaussian
processes with stationary increments have the same law when (and
only when) their spectral densities are equal a.e., we obtain that
\begin{eqnarray}
% \nonumber to remove numbering (before each equation)
   G_H &=& \{ C\in GL(n) : EB_H(t) B_H(s)^* = E(CB_H(t)) (CB_H(s))^*,\  s,t\in\bbR \} \nonumber  \\
   &=& \{ C\in GL(n) : (x^{-D}_{+}A + x^{-D}_{-}\overline{A}) (x^{-D}_{+}A + x^{-D}_{-}\overline{A})^* \nonumber \\
   & & \hspace{1.5in} = C (x^{-D}_{+}A + x^{-D}_{-}\overline{A}) (x^{-D}_{+}A +
   x^{-D}_{-}\overline{A})^* C^*,\ x\in\bbR\} \nonumber \\
   &=&  \{ C\in GL(n) : x^{-D} AA^* x^{-D^*} = C x^{-D} AA^* x^{-D^*} C^*,\ x>0 \} \nonumber \\
   &=& G_{H,1} \bigcap G_{H,2},
   \label{e:symmetry-group-OFBM-2}
\end{eqnarray}
where
\begin{eqnarray}
% \nonumber to remove numbering (before each equation)
  G_{H,1} &=& \{ C\in GL(n) : x^{-D} \Re(AA^*) x^{-D^*} = C x^{-D} \Re(AA^*) x^{-D^*} C^*,\ x>0 \} , \\
  G_{H,2} &=& \{ C\in GL(n) : x^{-D} \Im(AA^*) x^{-D^*} = C x^{-D} \Im(AA^*) x^{-D^*} C^*,\ x>0 \} .
\end{eqnarray}

Consider first the set $G_{H,1}$. Using the decomposition
(\ref{e:decomp-Re-AA*}) and working under the assumption
(\ref{e:full-rank}), we have that
\begin{eqnarray}
% \nonumber to remove numbering (before each equation)
  G_{H,1} &=&  \{ C\in GL(n) : x^{-D} S_R\Lambda_R^2 S_R^* x^{-D^*} = C x^{-D} S_R\Lambda_R^2 S_R^* x^{-D^*} C^*,\ x>0 \}
  \nonumber\\
   &=& \{ C\in GL(n) : (\Lambda_R^{-1} S_R^* x^D C x^{-D} S_R \Lambda_R) (\Lambda_R^{-1} S_R^* x^D C x^{-D} S_R \Lambda_R)^*
   = I,\ x>0 \}\nonumber \\
   &=& \{ C\in GL(n) : \Lambda_R^{-1} S_R^* x^D C x^{-D} S_R \Lambda_R\in O(n),\ x>0
   \}. \label{e:symmetry-group-OFBM-3}
\end{eqnarray}
Taking $x=1$ and using the fact that $S_R$ is orthogonal, $C\in
G_{H,1}$ necessarily has the form
\begin{equation}\label{e:C-form}
    C = S_R \Lambda_R S_R^* O S_R \Lambda_R^{-1} S_R^* = W O W^{-1}
\end{equation}
with $O\in O(n)$ (see also Remark \ref{r:on-O} below). Substituting
(\ref{e:C-form}) back into (\ref{e:symmetry-group-OFBM-3}), we can
now express $G_{H,1}$ as
\begin{eqnarray}
% \nonumber to remove numbering (before each equation)
  G_{H,1} &=& W
\{ O\in O(n) : \nonumber\\
& & \hspace{.3in} O (W^{-1} x^{-D} \Re(AA^*) x^{-D^*} W^{-1}) =
(W^{-1} x^{-D} \Re(AA^*) x^{-D^*} W^{-1}) O,\ x>0 \} W^{-1}
\nonumber\\
   &=& W \bigcap_{x > 0} G(\Pi_x) W^{-1}, \label{e:symmetry-group-OFBM-4}
\end{eqnarray}
where we use the definition (\ref{e:G-Pi}) of $G(\Pi_x)$, and
\begin{equation}\label{e:Pi_x}
    \Pi_x := W^{-1}  x^{-D} \Re(AA^*) x^{-D^*} W^{-1} = x^{-M}
    x^{-M^*}
\end{equation}
with
\begin{equation}\label{e:M}
    M = W^{-1} D W.
\end{equation}

\begin{remark}\label{r:on-O}
A simpler way to write (\ref{e:C-form}) and
(\ref{e:symmetry-group-OFBM-4}) would be to replace $W =
S_R\Lambda_R S_R^*$ by $S_R\Lambda_R$. Note that, with our choice,
$W$ is positive definite. The relation
(\ref{e:symmetry-group-OFBM-4}) then takes the form
(\ref{e:G_is_conjugate_to_compact_group}).
\end{remark}

The relation (\ref{e:symmetry-group-OFBM-4}) describes the first set
$G_{H,1}$ in the intersection (\ref{e:symmetry-group-OFBM-2}).
Instead of describing the second set $G_{H,2}$
separately, it is more convenient to think of the latter as imposing
additional conditions on the elements of $G_{H,1}$. In this regard,
observe first that, for any $y>0$,
\begin{equation}\label{e:G_H1-scaling}
    G_{H,1} = y^D G_{H,1} y^{-D},
\end{equation}
which simply follows by observing that the condition
$$
x^{-D} \Re(AA^*) x^{-D^*} = C x^{-D} \Re(AA^*) x^{-D^*} C^*,\ x>0,
$$
defining the set $G_{H,1}$, is equivalent to
$$
x^{-D} \Re(AA^*) x^{-D^*} = (y^D C y^{-D})x^{-D} \Re(AA^*) x^{-D^*}
(y^{-D^*} C^* y^{D^*}),\ x>0.
$$

Using the relation (\ref{e:G_H1-scaling}), $C\in G_{H,1}$
satisfies the relation
\begin{equation}\label{e:G_H2-define}
x^{-D} \Im(AA^*) x^{-D^*} = C x^{-D} \Im(AA^*) x^{-D^*} C^*,\ x>0,
\end{equation}
defining the set $G_{H,2}$, if and only if $C\in G_{H,1}$
satisfies the same relation (\ref{e:G_H2-define}) with $x=1$. Considering the form
(\ref{e:C-form}) of $C\in G_{H,1}$, this imposes additional
conditions on the orthogonal matrices $O$. Substituting
(\ref{e:C-form}) into the relation (\ref{e:G_H2-define}) with
$x=1$, we obtain that
$$
\Im(AA^*)  = W O W^{-1} \Im(AA^*) W^{-1} O^* W,
$$
i.e.,
$$
O W^{-1} \Im(AA^*) W^{-1} = W^{-1} \Im(AA^*) W^{-1} O,
$$
or
\begin{equation}\label{e:O-extra}
    O\in G(\Pi_I),
\end{equation}
where
\begin{equation}\label{e:Pi_I}
    \Pi_I = W^{-1} \Im(AA^*) W^{-1}.
\end{equation}
By the expressions (\ref{e:symmetry-group-OFBM-2}),
(\ref{e:symmetry-group-OFBM-4}) and the discussion above, we arrive at the following general result on the structure of symmetry groups
of OFBMs, and in particular, on the form of the conjugacy matrix
$W$.

\begin{theorem}\label{t:symmetry-group-OFBM}
Consider an OFBM given by the spectral representation
(\ref{e:spectral-repres-OFBM}), and suppose that the matrix $A$
satisfies the assumption (\ref{e:full-rank}). Then, its symmetry
group $G_H$ can be expressed as
\begin{equation}\label{e:symmetry-group-OFBM-main}
    G_H = W \Big( \bigcap_{x>0} G(\Pi_x) \cap G(\Pi_I) \Big) W^{-1},
\end{equation}
where $W$ is defined in (\ref{e:W}), and $\Pi_x$ and $\Pi_I$ are given
in (\ref{e:Pi_x}) and (\ref{e:Pi_I}), respectively.
\end{theorem}

The intersection over uncountably many $x>0$ in
(\ref{e:symmetry-group-OFBM-main}) can be replaced by a countable
intersection in a standard way. We have $O\in \cap_{x>0} G(\Pi_x)$
if and only if
\begin{equation}\label{e:G_H1-define}
    O x^{-M}x^{-M^*} = x^{-M}x^{-M^*} O,\ x>0.
\end{equation}
Writing $x^{-M} = \sum_{k=0}^\infty M^k (-\log(x))^k / k!$, the
relation (\ref{e:G_H1-define}) is equivalent to
$$
\sum_{k_1=0}^\infty \sum_{k_2=0}^\infty O M^{k_1} (M^*)^{k_2}
\frac{(-\log(x))^{k_1}}{k_1!} \frac{(-\log(x))^{k_2}}{k_2!} =
\sum_{k_1=0}^\infty \sum_{k_2=0}^\infty M^{k_1} (M^*)^{k_2} O
\frac{(-\log(x))^{k_1}}{k_1!} \frac{(-\log(x))^{k_2}}{k_2!}
$$
or, with $k_1 = k$, $k_1 + k_2 = m$,
$$
\sum_{m=0}^\infty O \sum_{k=0}^m M^k (M^*)^{m-k}
\frac{1}{k!(m-k)!} (-\log(x))^m = \sum_{m=0}^\infty \sum_{k=0}^m M^k
(M^*)^{m-k} O \frac{1}{k!(m-k)!} (-\log(x))^m.
$$
Equivalently,
\begin{equation}\label{e:O-Pi^m}
    O \Pi^{(m)} = \Pi^{(m)} O,\ m\geq 1,
\end{equation}
where
\begin{equation}\label{e:Pi^m}
    \Pi^{(m)} = \sum_{k=0}^m {m \choose k} M^k (M^*)^{m-k}.
\end{equation}
Theorem \ref{t:symmetry-group-OFBM} can now be reformulated as
follows.

\begin{theorem}\label{t:symmetry-group-OFBM-2}
Consider an OFBM given by the spectral representation
(\ref{e:spectral-repres-OFBM}), and suppose that the matrix $A$
satisfies the assumption (\ref{e:full-rank}). Then, its symmetry
group $G_H$ can be expressed as
\begin{equation}\label{e:symmetry-group-OFBM-main-2}
    G_H = W \Big( \bigcap_{m=1}^\infty G(\Pi^{(m)}) \cap G(\Pi_I) \Big) W^{-1},
\end{equation}
where $W$ is defined in (\ref{e:W}), and $\Pi^{(m)}$ and $\Pi_I$ are
given in (\ref{e:Pi^m}) and (\ref{e:Pi_I}), respectively.
\end{theorem}

\begin{remark}
Note that the matrix $\Pi_x$ in (\ref{e:Pi_x}) is positive definite. On the
other hand, the matrix $\Pi^{(m)}$ in (\ref{e:Pi^m}) is symmetric because so
are the terms
$$
{m\choose k} M^k (M^*)^{m-k} + {m\choose m-k} M^{m-k} (M^*)^{k}
$$
defining $\Pi^{(m)}$. However, $\Pi^{(m)}$ is not positive definite in
general. For example, with $\Re(AA^*)=I$ and normal $D$, we have
\begin{equation}\label{e:normal-D}
    \Pi^{(m)} = (D+D^*)^m,
\end{equation}
which is not positive definite (not even for $m=1$). Note
also that $\Pi_I$ is skew-symmetric, hence normal and diagonalizable.
\end{remark}

\section{On maximal and minimal symmetry groups}
\label{s:min-max}

An operator self-similar process $X$ is said to be of maximal
type, or elliptically symmetric, if its symmetry group $G_X$ is
conjugate to $O(n)$. At the other extreme, a zero mean (Gaussian)
o.s.s.\ process is said to be of minimal type if its symmetry
group is $\{I,-I\}$. We shall examine here these symmetry structures in the
case of OFBMs. First, we characterize maximal symmetry in terms of the spectral
parametrization of OFBMs. Second, we analyze minimal symmetry OFBMs
through a topological lens.

\subsection{OFBMs of maximal type}\label{s:maximal}

The following theorem is the main result of this subsection. Recall
the definition of single parameter OFBMs in Example
\ref{ex:single-param_OFBM}.

\begin{theorem}\label{t:maximal}
Consider an OFBM given by the spectral representation
(\ref{e:spectral-repres-OFBM}), and suppose that the matrix $A$
satisfies the assumption (\ref{e:full-rank}). If an OFBM is of
maximal type, then it is a single parameter OFBM up to a conjugacy
by a positive definite matrix. Moreover, this happens if and
only if
\begin{equation}\label{e:maximal-iff}
    \Im(AA^*) = 0,\quad -(D-dI) \Re(AA^*) = \Re(AA^*) (D^* - dI),
\end{equation}
for some real $d$.
\end{theorem}
\begin{remark}
Conversely, an OFBM which is a single parameter OFBM (up to a
positive definite conjugacy) is of maximal type (see Example
\ref{ex:single-param_OFBM}). We also point out that we have a proof of
the first claim in Theorem \ref{t:maximal} which does not make use
of spectral representations and dispenses with the assumption
(\ref{e:full-rank}). In the proof of Theorem \ref{t:maximal} we use spectral representations in order
to illustrate how the main results of Section \ref{s:groups} can
be used.
\end{remark}

\begin{proof} Since the OFBM is of maximal type, the representation $G_{H} = W_1 O(n) W^{-1}_1$ holds, where $W_1$ is positive definite. On the other hand, by (\ref{e:symmetry-group-OFBM-main}), we also have $G_{H} = W_2 {\mathcal O} W^{-1}_2$, where $W_2 $ is positive definite and ${\mathcal O}$ is a subgroup of the orthogonal group. Thus, by Lemma \ref{l:WO(n)W^(-1)=W2OW^(-1)2}, ${\mathcal O} = O(n)$. Therefore, in view of Proposition \ref{p:Centr=O(n)=>lambdaI}, maximal type occurs if and only if
\begin{equation}\label{e:min-max-Pi}
    \Pi_x = \lambda_x I,\ x>0, \quad \Pi_I = \lambda I, \quad \lambda_x,\lambda \in \bbR.
\end{equation}
Note that
\begin{equation}\label{e:min-max-Pi-I}
\Pi_I = \lambda I \quad \Leftrightarrow \quad \Im(AA^*) = \lambda
\Re(AA^*) \quad \Leftrightarrow \quad  \lambda = 0 \quad
\Leftrightarrow \quad \Im(AA^*) = 0.
\end{equation}
Moreover,
\begin{equation}\label{e:min-max-Pi-x}
\Pi_x = \lambda_x I,\ x>0 \quad \Leftrightarrow \quad \lambda_x
\Re(AA^*) = x^{-D} \Re(AA^*) x^{-D^*},
\end{equation}
which implies that, for any $x_1,x_2>0$,
$$
\lambda_{x_1x_2} \Re(AA^*) = (x_1x_2)^{-D} \Re(AA^*)
(x_1x_2)^{-D^*} = x_2^{-D} \lambda_{x_1} \Re(AA^*) x_2^{-D^*} =
\lambda_{x_1}\lambda_{x_2} \Re(AA^*).
$$
%Furthermore, for $\{x_k\}_{k} \subseteq (0,\infty)$ such that $x_k \rightarrow x \in (0,\infty)$,
%$$
%\lambda_{x_k}\Re(AA^*) = x^{-D}_{k}\Re(AA^*)x^{-D^*}_k \rightarrow x^{-D}\Re(AA^*)x^{-D^*} = \lambda_x \Re(AA^*).
%$$
Hence, under assumption (\ref{e:full-rank}), $\lambda_{x_1x_2} =
\lambda_{x_1}\lambda_{x_2}$, $x_1,x_2>0$. Moreover, the function $\log(\lambda_{\exp(\cdot)})$ is additive over $\bbR$, and it is measurable (since
it is continuous). As a consequence, by Theorem 1.1.8 in Bingham et al.\ \cite{bingham:goldie:teugels:1987}, p.\ 5, there exists a real $d$ such that
 $\log(\lambda_{\exp(\cdot)}) = -2d (\cdot)$, i.e., $\lambda_x = x^{-2d}$. In particular,
\begin{equation}\label{e:min-max-Pi-x-d}
x^{-D} \Re(AA^*) x^{-D^*} = x^{-2d} \Re(AA^*).
\end{equation}
Relations (\ref{e:min-max-Pi-I}) and (\ref{e:min-max-Pi-x-d}) imply
that the covariance structure of OFBM can be written as
$$
EB_H(t)B_H(s)^* = \int_\bbR \frac{e^{itx}-1}{ix}
\frac{e^{-isx}-1}{-ix} |x|^{-2d} \Re(AA^*) dx.
$$
In view of the relations (\ref{e:decomp-Re-AA*}) and (\ref{e:spectral-repres-single-OFBM}), this shows that $B_H$ is a
single parameter OFBM up to a conjugacy.

Finally, note from the above that $\Pi_x = \lambda_x I$, $x>0$, is
equivalent to $x^{-D} \Re(AA^*) x^{-D^*} = x^{-2d} \Re(AA^*)$ or
$x^{D-dI}\Re(AA^*) = \Re(AA^*) x^{-(D^*-dI)}$ for $x>0$ and some
real $d$. The latter is equivalent to $(D-dI)\Re(AA^*) = -
\Re(AA^*) (D^*-dI)$ for some real $d$. $\Box$
\end{proof}
\begin{remark}
Theorem 6 in Hudson and Mason \cite{hudson:mason:1982} shows that every maximal symmetry
o.s.s.\ process has an exponent of the form $h \hspace{0.5mm}I$,
$h \in \bbR$. For the case of OFBMs, the proof of Theorem
\ref{t:maximal} retrieves this result (see expression (\ref{e:min-max-Pi-x-d})).
Moreover, it is clear that, for a maximal symmetry OFBM $B_H$ (or, as a matter of fact,
for any maximal symmetry o.s.s.\ process) and for any $H \in {\mathcal E}(B_H)$,
$W^{-1}HW$ is normal, since $W^{-1}(H - hI)W \in so(n)$ (see also Section \ref{s:2-3} for further
results on the structure of exponents for dimensions $n=2$ and
$n=3$).
\end{remark}

\subsection{OFBMs of minimal type: the topologically general case}\label{s:min}

In view of Theorem \ref{t:symmetry-group-OFBM}, an OFBM is of minimal type if and only if
$$
\bigcap_{x > 0}G(\Pi_x) \cap G(\Pi_{I}) = \{I,-I\},
$$
and, in particular, if
\begin{equation}\label{e:centralizer_is_minimal}
\bigcap_{x > 0}G(\Pi_x) = \{I,-I\}.
\end{equation}
We shall focus here on the relation
(\ref{e:centralizer_is_minimal}) with the following related goals
in mind.

The first goal is to provide (practical) conditions for
(\ref{e:centralizer_is_minimal}) to hold and, hence, for an OFBM
to be of minimal type. This is a non-trivial problem.
The structure of $G(\Pi_x) $ depends on both the eigenvalues and the
invariant subspaces of $\Pi_x$, which are arbitrary in principle. Moreover,
their explicit calculation becomes increasingly difficult with
dimension. To shed light on (\ref{e:centralizer_is_minimal}), we
take up an idea from Lie group theory: a lot of information about
$M$ in the expression $\Pi_x = x^{-M}x^{-M^*}$ (see
(\ref{e:Pi_x})) is available at the vicinity of the identity in
the Lie group, i.e., as $x \rightarrow 1$. The general approach we take is
to study the behavior of the logarithm of $\Pi_x$ through the
Baker-Campbell-Hausdorff formula, i.e., by looking at the
associated Lie algebra. The characterization of the
behavior of the eigenvectors of $\Pi_x$ will then be retrieved by
turning back to the Lie group through the exponential map.

Initially, our conditions for the relation
(\ref{e:centralizer_is_minimal}) to hold are in terms of the
matrix $M$, and not directly in terms of $H$ and $A$. Our second
goal in this section is to show that these conditions on $M$ yield ``most" OFBMs
in terms of the parametrization $M$, and then relate them back to
$H$ and $A$. The term ``most" is in the topological sense, i.e.,
except on a meager set. This result should not be surprising: if
$\cap_{x>0}G(\Pi_x)$ has non-trivial structure, then this imposes
extra conditions on $M$ (or $D$, $W$) as in Section
\ref{s:maximal}. Though not surprising, formalizing this fact is
not straightforward, as shown here. This second goal leads to the main result of this section,
which, for the sake of clarity, we now briefly describe. In analogy with the assumption (\ref{e:eigen-assumption}), consider the set
\begin{equation}\label{e:mathcal_D}
{\mathcal D} = \Big\{D \in M(n,\bbR): -\frac{1}{2}< \Re(d_k) <
\frac{1}{2}, k=1,\hdots,n  \Big\},
\end{equation}
where $d_1,\hdots,d_n$ denote the charateristic roots of $D$. Theorem \ref{t:minimal_is_topol_general} below states the existence of a set ${\mathcal M} \subseteq M(n,\bbR)$ such that, for all $D$ and positive definite $W$ such that $W^{-1}DW \in {\mathcal M} \cap {\mathcal D}$, the OFBM with spectral parametrization $D$ and $\Re(AA^*):=W^2$ has minimal symmetry. Moreover, ${\mathcal M} \cap {\mathcal D}$ is an open set (of parameters), and it is dense in ${\mathcal D}$. Therefore, ${\mathcal M}^c \cap {\mathcal D}$ is a meager set. Conversely, every $M \in {\mathcal M} \cap {\mathcal D}$ gives a minimal symmetry OFBM through an appropriate spectral parametrization. In order to provide easy access to the mathematical content of this section, Remark \ref{r:heuristic_proof_identifiability} below contains a short heuristic proof of a weaker version of this claim, i.e., that OFBMs have identifiable (unique) exponents for every parametrization except on a meager set.

The rest of this section is dedicated to developing these ideas, as well as the framework behind them. Hereinafter, unless otherwise stated, we impose no restrictions on the eigenvalues of $M$, i.e., the expression
$\Pi_x = x^{-M}x^{-M^*}$ is taken for any $M \in M(n,\bbR)$.
Consider the decomposition of the latter space into the direct sum
\begin{equation}\label{e:decomp_M}
M(n,\bbR) = {\mathcal S}(n,\bbR) \oplus so(n).
\end{equation}
Correspondingly, for $M \in M(n,\bbR)$, denote
\begin{equation}\label{e:M=S+L}
M = S + L,
\end{equation}
where $S = (M+M^*)/2$, $L=(M-M^*)/2$ are, respectively, the symmetric and skew-symmetric parts
of $M$.

The next proposition shows that for an appropriately chosen $M$,
the centralizer of the family $\Pi_x$ is minimal. In the proof, the symbol $[\cdot,\cdot]$ denotes
the commutator. Since the point $x=1$ is a singularity in the
sense that all the information about $M$ from $\Pi_x = x^{-M}x^{-M^*}$
 is lost at it, the idea
is to analyze the behavior of $\Pi_x$ for $x$ in a close vicinity
of 1. The proof of Proposition \ref{p:invariant_subspaces_of_Pi(x)_change} is based on the idea that, for a matrix parameter $M = S+L$, $S \in {\mathcal S}_{\neq} \subseteq {\mathcal S}(n,\bbR)$  and $L \in so(n)$, the existence of non-trivial solutions for the matrix equations $O \Pi_x = \Pi_x O$, $x > 0$, implies the coincidence of some invariant subspace of $S$ and $L$, which, as we will see afterwards, is a very special situation in the topological sense.

\begin{proposition}\label{p:invariant_subspaces_of_Pi(x)_change}
Let $M = S + L$, where $S \in {\mathcal S}_{\neq}$ (see expression (\ref{e:Sneq})) and $L \in so(n)$ do not share $k$-dimensional real invariant subspaces, $k = 1, \hdots, n-1$. Then, $\bigcap_{x > 0}G(\Pi_x) = \{I,-I\}$.
\end{proposition}
\begin{proof}
Note that $M+M^* = 2S$, and that
$$
[M,M^*]= MM^*-M^*M = (S+L)(S-L) - (S-L)(S+L) = 2[L,S].
$$
Since the mapping $M \mapsto \exp(M)$ is a $C^{\infty}$
homeomorphism of some neighborhood of 0 in the Lie algebra of
$GL(n,\bbR)$ onto some neighborhood $U$ of the identity $I$ in
$GL(n,\bbR)$, then its inverse function $\textnormal{Log}$ is
well-defined on $U$. Therefore, by the Baker-Campbell-Hausdorff
formula, for small enough $\log(x)$ we have
$$
\textnormal{Log}(\exp(-\log(x)M)\exp(-\log(x)M^*)) =
-\log(x)(M+M^*) + \frac{1}{2}[-\log(x)M,-\log(x)M^*]
$$
$$
 +
O(\log^3(x)) = -\log(x)(M+M^*) + \log^2(x) \frac{1}{2}[M,M^*] + O(\log^3(x))
$$
(see Hausner and Schwartz \cite{hausner:schwartz:1968}, p.\ 63 and pp.\ 68-69).

We would like to show that there exists $\delta > 0$ such that the symmetric matrices
\begin{equation}\label{e:2nd_order}
-\frac{1}{\log(x)}\textnormal{Log}(x^{-M}x^{-M^*}) = (M+M^*) - \frac{\log(x)}{2}  [M,M^*] + O(\log^2(x)), \quad x \in B(1,\delta) \backslash\{1\},
\end{equation}
do not share $k$-dimensional real invariant subspaces with the symmetric matrix $M+M^*$, $k=1,\hdots,n-1$.

Without loss of generality, assume that there exists $\{x_i\}$, $x_i \rightarrow 1$, such that each associated expression in (\ref{e:2nd_order}) shares a $k(i)$-dimensional invariant subspace with $M+M^*$. Since the number of possible real invariant subspaces of $M+M^* \in {\mathcal S}_{\neq}(n,\bbR) $ is finite (by Corollary \ref{c:invar_subs=>eigenvec}, they are generated by real eigenvectors of $M+M^*$), by passing to a subsequence if necessary, we can assume that each
$$
-\frac{1}{\log(x_i)} \textnormal{Log}(x^{-M}_{i}x^{-M^*}_{i})
$$
shares the same invariant subspace with $M+M^*$. Write
$$
M+M^* = O(2\Upsilon)O^*,
$$
where $\Upsilon$ is diagonal and assume without loss of generality that the invariant subspace in question is $\textnormal{span}_{\bbR}\{o_1,\hdots,o_k\}$ (i.e., the first $k$ columns of $O \in O(n)$). This implies that we can write
$$
O(2 \Upsilon)O^* - \frac{\log(x_i)}{2}((O\Upsilon O^*)L - L(O\Upsilon O^*)) + O(\log^2(x_i)) = OJ_iO^*,
$$
where by Corollary \ref{c:invar=>zeroes_in_S_and_L} $J_i$ is \textit{block}-diagonal of the form
$$
J_i = \left(\begin{array}{cc}
J^{i}_{11} & {\mathbf 0} \\
{\mathbf 0} & J^{i}_{22} \\
\end{array}\right), \quad J^{i}_{11} \in M(k,\bbR), \quad J^{i}_{22} \in M(n-k,\bbR).
$$
Therefore,
$$
2 \Upsilon - \frac{\log(x_i)}{2}(\Upsilon (O^*LO)- (O^*LO)\Upsilon) = J_i + O(\log^2(x_i)).
$$
If $L_2 := O^*LO$, then
\begin{equation}\label{e:D(O^*LO)- (O^*LO)D}
\Upsilon L_2- L_2 \Upsilon = (J_i - 2 \Upsilon) \Big(-\frac{2}{\log(x_i)} \Big) + O(\log(x_i)).
\end{equation}
Note that the term $J_i - 2 \Upsilon$ is still block-diagonal of the same form as $J_i$. For notational simplicity, we still write $J_i - 2 \Upsilon = \textnormal{diag}(J^{i}_{11},J^{i}_{22})$. Therefore, the right-hand side of (\ref{e:D(O^*LO)- (O^*LO)D}) must converge as $x_i \rightarrow 1$. Denote the limit by
\begin{equation}\label{e:J_blockdiag}
J = \left(\begin{array}{cc}
J_{11} & {\mathbf 0} \\
{\mathbf 0} & J_{22} \\
\end{array}\right),
\end{equation}
where $J_{11} \in {\mathcal S}(k,\bbR)$, $J_{22} \in {\mathcal S}(n-k,\bbR)$.
Denote $L_2 = (l_{ij})_{i,j=1,\hdots,n}$, $\Upsilon = \textnormal{diag}(\lambda_1,\hdots,\lambda_n) $. Then
\begin{equation}\label{e:D L_2- L_2 D}
\Upsilon L_2- L_2 \Upsilon
= \left(\begin{array}{cccccc}
0 & l_{12}(\lambda_1 - \lambda_2) & l_{13}(\lambda_1 - \lambda_3) & l_{14}(\lambda_1 - \lambda_4) & \hdots & l_{1n}(\lambda_1 - \lambda_n) \\
  &    0    & l_{23}(\lambda_2 - \lambda_3) & l_{24}(\lambda_2 - \lambda_4) & \hdots & l_{2n}(\lambda_2 - \lambda_n) \\
  &        & 0 & l_{34}(\lambda_3 - \lambda_4) & \hdots & l_{3n}(\lambda_3 - \lambda_n) \\
  &        &  & \ddots & \hdots & \vdots \\
  &        &  &         &  &  l_{n-1,n}(\lambda_{n-1} - \lambda_n) \\
  &        &  &          &   & 0 \\
\end{array}\right),
\end{equation}
where the entries below the main diagonal are equal to the corresponding ones above the main diagonal. From expressions (\ref{e:D L_2- L_2 D}) and (\ref{e:J_blockdiag}), we conclude that the upper right (non-square) block is zero, i.e.,
$$
\left(\begin{array}{cccccc}
l_{1,k+1}(\lambda_1 - \lambda_{k+1}) & l_{1,k+2}(\lambda_1 - \lambda_{k+2}) & \hdots & l_{1,n}(\lambda_1 - \lambda_{n}) \\
l_{2,k+1}(\lambda_2 - \lambda_{k+1}) & l_{2,k+2}(\lambda_2 - \lambda_{k+2}) & \hdots & l_{2,n}(\lambda_2 - \lambda_{n}) \\
\vdots & \vdots & \hdots &  \vdots  \\
l_{k,k+1}(\lambda_k - \lambda_{k+1}) & l_{k,k+2}(\lambda_k - \lambda_{k+2}) & \hdots & l_{k,n}(\lambda_k - \lambda_{n}) \\
\end{array}\right) = {\mathbf 0}.
$$
Therefore, since the $\lambda$'s are pairwise different,
$$
L_{12}:= \left(\begin{array}{cccccc}
l_{1,k+1}  & l_{1,k+2}  & \hdots & l_{1,n} \\
l_{2,k+1}  & l_{2,k+2}  & \hdots & l_{2,n} \\
\vdots & \vdots & \hdots &  \vdots  \\
l_{k,k+1}  & l_{k,k+2}  & \hdots & l_{k,n} \\
\end{array}\right) = {\mathbf 0}.
$$
Therefore,
$$
L = O
\left(\begin{array}{cc}
L_{11} & L_{12} \\
L_{21} & L_{22}
\end{array}\right)O^* =
O \left(\begin{array}{cc}
L_{11} & {\mathbf 0} \\
{\mathbf 0} & L_{22}
\end{array}\right)O^*
$$
(note that $L_{11}$ and $L_{22}$ may contain zeroes off the main diagonal but this is inconsequential). Therefore, $\textnormal{span}_{\bbR}\{o_1,\hdots,o_k\}$ is also an invariant subspace of $L$ (contradiction).

So, fix $x_0 \in B(1,\delta)\backslash\{1\}$ such that $\textnormal{Log}(x^{-M}_{0}x^{-M^*}_{0}) = -\log(x_0) (M+M^*) + \frac{1}{2}\log^2(x_0)[M,M^*] + O(\log^3(x_0))$ has no $k$-dimensional real invariant subspaces in common with $M+M^*$ or, equivalently, $-\log(x_0) (M+M^*)$, $k = 1,\hdots,n-1$. Moreover, by Lemma \ref{l:conv_of_eigenvalues}, we can assume that, for $x \in B(1,\delta)\backslash\{1\}$, $\textnormal{Log}(x^{-M}x^{-M^*}) \in {\mathcal S}_{\neq}(n,\bbR)$.

% Then, since the eigenvectors of $\textnormal{Log}(x^{-M}x^{-M^*})$ must converge (in the sense of Lemma \ref{l:Pk_conv_to_P}) to those of $M+M^*$ %(or, equivalently, $-\log(x) (M+M^*)$), then

We claim that there exists $x_1 \in B(1,\delta)\backslash\{1,x_0\}$ such that $\textnormal{Log}(x^{-M}_1 x^{-M^*}_1)$ shares no $k$-dimensional real invariant subspaces with $\textnormal{Log}(x^{-M}_0 x^{-M^*}_0)$. This can be proved by an argument similar to the one by contradiction used above in this proof. In fact, assume that there exists a sequence $\{x_{i}\}$ such that $\textnormal{Log}(x^{-M}_i x^{-M^*}_i )$ shares some $k(i)$-dimensional real invariant subspaces with $\textnormal{Log}(x^{-M}_0 x^{-M^*}_0 )$. Since the latter is symmetric, then by Corollary \ref{c:invar_subs=>eigenvec} its real invariant subspaces are generated by $k$ eigenvectors of $\textnormal{Log}(x^{-M}_0 x^{-M^*}_0 )$. Denote a basis of orthonormal eigenvectors of the latter by $v_1,\hdots,v_{n}$. Since $\textnormal{Log}(x^{-M}_0 x^{-M^*}_0) \in {\mathcal S}_{\neq}(n,\bbR)$ and the number of possible subspaces of the form $\textnormal{span}_{\bbR}\{v_{j_1},\hdots,v_{j_k}\}$ is finite, then we can assume that the shared invariant subspace is the same for all $i$. For notational simplicity, write it as $\textnormal{span}_{\bbR}\{v_{1},\hdots,v_{k}\}$; complete this basis with orthonormal vectors $v_{k+1},\hdots,v_{n}$ and let $P = (v_{1},\hdots,v_{k},v_{k+1},\hdots,v_{n}) \in O(n)$. Then, by Corollary \ref{c:invar=>zeroes_in_S_and_L}
$$
(M+M^*) + O(\log(x_i)) = - \frac{1}{\log(x_i)}\textnormal{Log}(x^{-M}_i x^{-M^*}_i ) = P \textnormal{diag}(J^{i}_{11},J^{i}_{22})P^*.
$$
Thus, $\lim_{i \rightarrow \infty} J^{i}_{11}$, $\lim_{i \rightarrow \infty} J^{i}_{22}$ exist, which implies that $\textnormal{span}_{\bbR}\{v_{1},\hdots,v_{k}\}$ is a real invariant subspace of $M+M^*$ (contradiction).

Consequently, we also have that $x^{-M}_1 x^{-M^*}_1$ and $x^{-M}_0 x^{-M^*}_0$ share no $k$-dimensional real invariant subspaces (since the eigenvectors are the same as those of $\textnormal{Log}(x^{-M}_1 x^{-M^*}_1)$ and $\textnormal{Log}(x^{-M}_0 x^{-M^*}_0)$, respectively, and by Corollary \ref{c:invar_subs=>eigenvec}, the $k$-dimensional real invariant subspaces of symmetric matrices are each generated by a set of $k$ real eigenvectors). Thus, by Lemma \ref{l:irreducible}, the conclusion follows, i.e., the only orthogonal matrices that commute with both $x^{-M}_0 x^{-M^*}_0$ and $x^{-M}_1 x^{-M^*}_1$ are $\pm I$.
$\Box$\\
\end{proof}

Bearing in mind Proposition \ref{p:invariant_subspaces_of_Pi(x)_change}, we would like to construct a set of matrices $M$ based on which we can apply the proposition, and whose topology we can characterize. We now take a closer look at an appropriate set of skew-symmetric matrices.

\begin{definition}\label{def:L_invar}
Let $o_1,\hdots, o_n$ be a real orthonormal basis of $\bbR^n$. Let
$$
{\mathcal L}_{\textnormal{invar}}(o_1,\hdots, o_n) = \{L \in so(n):
$$
\begin{equation}\label{e:L_invar}
\textnormal{there exists a subset $o_{j_1},\hdots,o_{j_k}$, $k < n$, that generates a real invariant subspace of $L$}\}.
\end{equation}
\end{definition}

\begin{example}\label{ex:example_L_invar}
If $a,b,c,d \in \bbR$, then
$$
L := \left(\begin{array}{ccccc}
0 & a & b & & \\
-a & 0 & c & & \\
-b & -c & 0 & & \\
   &    &   & 0 & d \\
   &    &   & -d & 0 \\
\end{array}\right) \in {\mathcal L}_{\textnormal{invar}}(e_1,e_2,e_3,e_4,e_5)
$$
since the real subspace $\textnormal{span}_{\bbR}\{e_1,e_2,e_3\}$ is invariant with respect to $L$ (the same being true for $\textnormal{span}_{\bbR}\{e_4,e_5\}$).
\end{example}

\begin{remark}
The representations of the sets ${\mathcal L}_{\textnormal{invar}}$ may use different arguments (vectors). For instance, ${\mathcal L}_{\textnormal{invar}}(v_1,v_2) = \{\mathbf{0}\}$ for every orthonormal pair $v_1, v_2 \in \bbR^2$. However, this is inconsequential for the developments in this section.
\end{remark}

Proposition \ref{p:open-dense_direct_sum} below establishes the
topological properties of the class of ``well-behaved" exponents
$M$, i.e., those that will eventually be associated with minimal
symmetry OFBMs. Its proof is based on the next two lemmas.

\begin{lemma}\label{e:open-dense_separate}
\begin{enumerate}
\item [(i)] The set ${\mathcal S}_{\neq}$ is an open, dense set in
(the relative topology of) ${\mathcal S}(n,\bbR)$.

\item [(ii)] For any orthonormal vectors
$o_1,\hdots,o_n$ in $\bbR^n$, the set $({\mathcal L}_{\textnormal{invar}}(o_1,\hdots,o_n))^{c}$ is an open, dense set in (the relative topology of) $so(n)$.
\end{enumerate}
\end{lemma}
\begin{proof}
$(i)$ Openness stems from the fact that, for $\{S_{k}\} \subseteq {\mathcal S}(n,\bbR)$ such that $S_k \rightarrow S_0$, by Lemma \ref{l:conv_of_eigenvalues}, the eigenvalues of $S_k$ converge to those of $S_0$. Thus, for large enough $k$, the latter are pairwise distinct.

For denseness, take $S_0 \in ({\mathcal S}_{\neq})^c$. One can obtain a sequence $\{S_{k}\} \subseteq {\mathcal S}_{\neq}$ such that $S_k \rightarrow S_0$ simply by appropriately perturbing the eigenvalues of $S_0$, for a fixed conjugacy $O \in O(n)$ of eigenvectors of $S_0$.

$(ii)$ By contradiction, fix a real orthonormal basis $o_1,\hdots, o_n$ and assume that $({\mathcal L}_{\textnormal{invar}})^c:=({\mathcal L}_{\textnormal{invar}}(o_1,\hdots,o_n))^c$ is not open. Then there exists $L \in ({\mathcal L}_{\textnormal{invar}})^c$ such that, for some $\{L_i\}_{i \in \bbN} \subseteq {\mathcal L}_{\textnormal{invar}}$, $L_i \rightarrow L$. Since there are finitely many $k$-tuples $(o_{j_1},\hdots,o_{j_k})$, $k = 1,\hdots,n-1$, then we can extract a (convergent) subsequence $\{L_{i'}\}$ for which all $L_{i'}$'s share one invariant subspace $\textnormal{span}_{\bbR}\{o_{j_1},\hdots,o_{j_k}\}$ (i.e., $k$ is not a function of $i$). This means that we can form an orthogonal matrix $O := (o_{j_1},\hdots,o_{j_k}, o_{j_{k+1}},\hdots,o_{j_{n}})$, where $\textnormal{span}_{\bbR}\{o_{j_{k+1}},\hdots,o_{j_{n}}\} = (\textnormal{span}_{\bbR}\{o_{j_1},\hdots,o_{j_k}\})^{\perp}$, and by Corollary \ref{l:invar=>zeroes} write
$$
L_{i} = O \left(\begin{array}{cc}
L^{i}_{11} & {\mathbf 0} \\
{\mathbf 0}  & L^{i}_{22}
\end{array}\right) O^*,
$$
where $L^{i}_{11} \in so(k), L^{i}_{22} \in so(n-k)$. Since $\lim_{i \rightarrow \infty} L^{i}_{11}, \lim_{i \rightarrow \infty} L^{i}_{22}$ must exist, then $L \in {\mathcal L}_{\textnormal{invar}}$ (contradiction).

As for denseness, once again fix a real orthonormal basis $o_1,\hdots, o_n$. Take any $L \in {\mathcal L}_{\textnormal{invar}}:={\mathcal L}_{\textnormal{invar}}(o_1,\hdots,o_n)$ and write it in a block-diagonal form as $L = O \textnormal{diag}(L_{11},\hdots,L_{jj})O^*$, where $L_{11},\hdots,L_{jj}$ are skew-symmetric matrices. Now form the sequence of matrices $\{L_i\}$ by replacing all the zero entries above the main diagonal of $L$ with $1/i$, and correspondingly, the zero entries below the main diagonal with $-1/i$ (this may include entries in the blocks $L_{11},\hdots,L_{jj}$ themselves). Then, by Corollary \ref{c:invar=>zeroes_in_S_and_L}, we must have $L_{i} \in {\mathcal L}^c_{\textnormal{invar}}$, and $L_i \rightarrow L$.
$\Box$\\
\end{proof}

We now define a correspondence (set-valued function) that maps the
set ${\mathcal S}_{\neq}$ into sets of skew-symmetric matrices in the classes (\ref{e:L_invar}).
\begin{definition}\label{def:L}
Let ${\mathcal P}$ be the class of all subsets of a set. Define the correspondence (set-valued function)
$$
l: {\mathcal S}_{\neq} \rightarrow {\mathcal P}(so(n)),
$$
\begin{equation}\label{e:S_mapsto_l(S)}
S \mapsto l(S) =
({\mathcal L}_{\textnormal{invar}}(o_1,\hdots,o_n))^c
\end{equation}
where $o_1,\hdots,o_n$ represent orthonormal eigenvectors of $S$.
\end{definition}

\begin{remark}
The correspondence $l(\cdot)$ is well-defined. In fact, for $S \in {\mathcal S}_{\neq}(n,\bbR)$ and an associated basis of eigenvectors $o_1,\hdots,o_{n}$, all the possible representations of ${\mathcal L}_{\textnormal{invar}}(o_1,\hdots,o_n)$ are of the form ${\mathcal L}_{\textnormal{invar}}(\pm o_1,\hdots,\pm o_n)$. Now note that, for any subset $o_{j_1},\hdots, o_{j_k}$, $1 \leq k < n $, $\textnormal{span}_{\bbR}\{ o_{j_1},\hdots,o_{j_k} \}=\textnormal{span}_{\bbR}\{\pm o_{j_1},\hdots, \pm o_{j_k}\}$, i.e., invariant subspaces can equivalently be expressed in either basis.
\end{remark}

In the next lemma, we show that the graph of the correspondence (\ref{e:S_mapsto_l(S)}) is open and dense in the set $({\mathcal S}_{\neq},so(n))$. The topology under consideration is that generated by the open rectangles
\begin{equation}\label{e:rectangles}
B_{{\mathcal S}(n,\bbR)}(S_0,\varepsilon_{1}) \times B_{so(n)}(L_0,\varepsilon_{2}), \quad \varepsilon_1, \varepsilon_2 > 0,
\end{equation}
where
$$
B_{{\mathcal S}(n,\bbR)}(S_0,\varepsilon_{1}):= \{M \in M(n,\bbR):\norm{M - S_0}< \varepsilon_{1}\} \cap {\mathcal S}(n,\bbR),
$$
$$
B_{so(n)}(L_0,\varepsilon_{2}):= \{M \in M(n,\bbR):\norm{M - L_0}< \varepsilon_{2}\} \cap so(n)
$$
and $\norm{\cdot}$ is the spectral norm.

\begin{lemma}\label{l:graph(L)_is_open}
Let $l(\cdot)$ be the correspondence in Definition \ref{def:L}. Then $\textnormal{Graph}(l):= \{(S,L): S \in {\mathcal S}_{\neq}, L \in l(S)\}$ is open and dense in $({\mathcal S}(n,\bbR),so(n))$.
\end{lemma}
\begin{proof}
Openness is a consequence of the fact that, if $S_0 \in {\mathcal S}_{\neq}$, then, as $S_i \rightarrow S_0$, the eigenvalues of $S_i$ converge to those of $S_0$ (in the sense of Lemma \ref{l:conv_of_eigenvalues}). Indeed, assume by contradiction that there exists $(S_0,L_0) \in \textnormal{Graph}(l)$ such that, for some sequence $(S_i,L_i) \notin \textnormal{Graph}(l)$,
$$
(S_i,L_i) \rightarrow (S_0,L_0).
$$
Note that there cannot be a subsequence $\{S_{i'}\} \subseteq {\mathcal S}^{c}_{\neq}$ such that $S_{i'} \rightarrow S_0$ (since this contradicts
the openness of ${\mathcal S}_{\neq}$ established in Lemma \ref{e:open-dense_separate}). Thus, we can assume that $\{S_i\}_{i \in \bbN} \subseteq {\mathcal S}_{\neq}$. Consequently, we must have $L_{i} \in {\mathcal L}_{\textnormal{invar}}(o^{i}_1,\hdots,o^{i}_n)$, where $o^{i}_1,\hdots,o^{i}_n$ is a basis of real eigenvectors of $S_i$. Since $S_0 \in {\mathcal S}_{\neq}$, then by Lemma \ref{l:Pk_conv_to_P}, we can assume that $o^{i}_1 \rightarrow o_1 ,\hdots,o^{i}_n  \rightarrow o_n$, where $o_1,\hdots,o_n$ is a basis of real eigenvectors of $S_0$. Therefore, we can write
$$
L_i = (o^{i}_1,\hdots,o^{i}_n)K_i (o^{i}_1,\hdots,o^{i}_n)^{*}, \quad i \in \bbN,
$$
where, by Definition \ref{ex:example_L_invar} and Corollary \ref{c:invar=>zeroes_in_S_and_L}, possibly by a permutation of the vectors $o^{i}_1,\hdots,o^{i}_n$ the matrix $K_i \in so(n)$ can be made block-diagonal (see also Example \ref{def:L_invar}). Define
$$
k^* = \min\{k = 1,\hdots,n-1: \textnormal{infinitely many $L_{i}$'s have a $k$-dimensional real invariant subspace} \},
$$
i.e., the minimal non-trivial dimension for invariant subspaces over infinitely many terms of the sequence $\{L_i\}$. Now for the associated sequence of vectors $o^{i}_1,\hdots, o^{i}_n$, for each $i$ choose one subset of indices $j_1(i),\hdots,j_{k^*}(i) \subseteq \{1,\hdots,n\}$ such that $o^{i}_{j_1(i)},\hdots,o^{i}_{j_{k^*}(i)}$ generates a real invariant subspace of $L_i$ (there may be more than one choice, but this is inconsequential). Since there are at most ${n \choose k^*}$ such choices, by passing to a subsequence if necessary, we can fix a subset of indices $j_1,\hdots,j_{k^*}$ such that $o^{i}_{j_1},\hdots,o^{i}_{j_{k^*}}$ generates a real invariant subspace of $L_i$ for every $i$ in this (sub)sequence. Since we can change at will the order of columns in the conjugacy matrix, without loss of generality we can assume that $j_1 = 1,\hdots,j_{k^*} = k^*$. Therefore,
$$
L_i = O_i \textnormal{diag}(L^{i}_{11},L^{i}_{22})O^*_i, \quad i \in \bbN,
$$
where $L^{i}_{11} \in so(k^*),L^{i}_{22} \in so(n-k^*)$. Since $O_i = (o^{i}_{1},\hdots,o^{i}_{n}) \rightarrow (o_{1},\hdots,o_{n})$ and $L_i \rightarrow L_0$, then the limits $\lim_{i \rightarrow \infty}L^{i}_{11}$, $\lim_{i \rightarrow \infty}L^{i}_{22}$ exist and thus $L_0 \in {\mathcal L}_{\textnormal{invar}}(o_1,\hdots,o_n)$. Therefore, $(S_0,L_0) \notin \textnormal{Graph}(l)$ (contradiction).

We now show denseness. Assume $(S_0,L_0) \notin \textnormal{Graph}(l) $. We will break up the argument into cases. Let $o_1,\hdots,o_n$ be a basis of real orthonormal eigenvectors of $S_0$.

\begin{itemize}
\item $S_0 \in {\mathcal S}_{\neq}$: Then $L_0 \in {\mathcal L}_{\textnormal{invar}}(o_1,\hdots,o_n)$. We can apply the same argument as in the proof of Lemma \ref{e:open-dense_separate} and obtain the sequence $(S_0,L_i) \in \textnormal{Graph}(l) $, $(S_0, L_i) \rightarrow (S_0,L_0)$.
\item $S_0 \notin {\mathcal S}_{\neq},L_0 \in ({\mathcal L}_{\textnormal{invar}}(o_1,\hdots,o_n))^c$: Then generate a sequence $S_i = O D_i O^*$ by appropriately perturbing the repeated eigenvalues of $S_0$ so that $S_i \in {\mathcal S}_{\neq}$ and $S_i \rightarrow S_0$. Thus, $(S_i, L_0) \in \textnormal{Graph}(l)$ (since $S_i$ and $S_0$ share the eigenvector basis $o_1,\hdots,o_n$) and $(S_i, L_0)\rightarrow (S_0,L_0)$.
    \item $S_0 \notin {\mathcal S}_{\neq},L_0 \in {\mathcal L}_{\textnormal{invar}}(o_1,\hdots,o_n)$: As in the previous case, generate a sequence $S_i = O D_i O^*$ by appropriately perturbing the repeated eigenvalues of $S_0$ so that $S_i \in {\mathcal S}_{\neq}$ and $S_i \rightarrow S_0$. Without loss of generality, assume that the subset of vectors $o_1,\hdots, o_k$ generates a real invariant subspace with respect to $L_0$. Now apply the same argument as in the proof of Lemma \ref{e:open-dense_separate} and obtain the sequence $(S_i,L_i) \in \textnormal{Graph}(l) $ (since $S_i$ and $S_0$ share the eigenvector basis $o_1,\hdots,o_n$), with $(S_i, L_i) \rightarrow (S_0,L_0)$.
\end{itemize}
$\Box$\\
\end{proof}

\begin{remark}
Regarding the last part of the proof of Lemma \ref{l:graph(L)_is_open} (on denseness), when $S_0 \notin {\mathcal S}_{\neq}(n,\bbR)$ there are infinitely many choices of bases of orthonormal eigenvectors of $S_0$. Thus, $L_0$ is in ${\mathcal L}_{\textnormal{invar}}(o_1,\hdots,o_n )$ or not depending on the particular basis $o_1,\hdots,o_n$ chosen.
\end{remark}

The next proposition puts $\textnormal{Graph}(l)$ back into the
original space $M(n,\bbR)$ in the form of a direct sum, rephrases
the topological statement of Lemma \ref{l:graph(L)_is_open}, and
connects the latter to the problem of proving
(\ref{e:centralizer_is_minimal}).
\begin{proposition}\label{p:open-dense_direct_sum}
Let
\begin{equation}\label{e:Mset}
{\mathcal M}=\{M \in M(n,\bbR): M \in {\mathcal S}_{\neq} \oplus
l({\mathcal S}_{\neq} )\}.
\end{equation}
Then,
\begin{itemize}
\item [(i)] ${\mathcal M}$ is an open, dense subset of
$M(n,\bbR)$. Consequently, ${\mathcal M}^c$ is a meager set and
${\mathcal M}$ is a $n^2$-dimensional $C^{\infty}$ manifold in
$\bbR^{n^2} \cong M(n,\bbR)$. \item [(ii)] relation
(\ref{e:centralizer_is_minimal}) holds for all $M \in {\mathcal
M}$.
\end{itemize}
\end{proposition}
\begin{proof}
We first show part $(i)$. Define the norm $\norm{M}_{\oplus} = \norm{S}+\norm{L}$, where $\norm{\cdot}$ is the spectral matrix norm. Expression (\ref{e:decomp_M}) implies that $\norm{\cdot}_{\oplus}$ is well-defined. By the equivalence of matrix norms, it suffices to show $(i)$ with respect to $\norm{\cdot}_{\oplus}$.

Assume by contradiction that ${\mathcal M}$ is not open. Then there is some $M_0 = S_0 + L_0 \in {\mathcal M}$ and a sequence $\{M_{k}\} \subseteq {\mathcal M}^{c}$ such that $\norm{M_k - M_0}_{\oplus} \rightarrow 0$. However, such convergence holds if and only if $\norm{S_k - S_0}\rightarrow 0$ and $\norm{L_k - L_0}\rightarrow 0$. Since, for each $k$, either $S_k \notin {\mathcal S}_{\neq}$ or $L_k \notin l(S_k)$, then this contradicts the openness of $\textnormal{Graph}(l)$ (Lemma \ref{l:graph(L)_is_open}). Denseness can be addressed in a similar fashion, and the geometric statement
is immediate.

Part $(ii)$ is a consequence of Proposition \ref{p:invariant_subspaces_of_Pi(x)_change}. $\Box$\\
\end{proof}

\begin{example}
To construct an example of $M \in {\mathcal M}$, we turn to
the case when $n=3$. Take $S =
\textnormal{diag}(s_1,s_2,s_3)$, where the (real) eigenvalues are
pairwise different. Then, take any $L \in so(3)\backslash\{0\}$
\textit{not} having a real invariant subspace generated by a 1 or 2-tuple of Euclidean canonical vectors $e_1$,
$e_2$, $e_3$. In other words, we
can\textit{not} take a matrix $L$ of one of the forms
$$
\left(\begin{array}{ccc}
0 & a & 0\\
-a & 0 & 0\\
0 & 0 & 0
\end{array}\right), \quad
\left(\begin{array}{ccc}
0 & 0 & b\\
0 & 0 & 0\\
-b & 0 & 0
\end{array}\right),\quad
\left(\begin{array}{ccc}
0 & 0 & 0\\
0 & 0 & c\\
0 & -c & 0
\end{array}\right),
$$
where $a,b,c \in \bbR$. Now set $M = S+L$.
\end{example}

In order to make the claim about the general minimality of the symmetry groups of OFBMs, we need to restrict the parameter space,
as in (\ref{e:eigen-assumption}). For this purpose, we consider the set ${\mathcal D}$ in (\ref{e:mathcal_D}). The following is the main result of this section. It shows that, except possibly when the parametrization is taken on a meager set,
OFBMs are of minimal symmetry.
\begin{theorem}\label{t:minimal_is_topol_general}
For all $D, W \in M(n,\bbR)$, where $W$ is positive definite, and such that
$$
W^{-1}DW \in {\mathcal M} \cap {\mathcal D},
$$
the associated OFBM with spectral parametrization $D$ and $\Re(AA^*):= W^2$ has minimal symmetry. The set ${\mathcal M} \cap {\mathcal D}$ is open, and,
in particular, it is an $n^2$-dimensional $C^{\infty}$ manifold in $\bbR^{n^2} \cong M(n,\bbR)$. Moreover, it is also a dense subset of ${\mathcal D}$
(in the relative topology of ${\mathcal D}$). As a consequence, ${\mathcal M}^{c} \cap {\mathcal D}$ is a meager set.

Conversely, every $M \in {\mathcal M} \cap {\mathcal D}$ gives rise to a minimal symmetry OFBM through the spectral parametrization $D:=M$, $W:=I$.
\end{theorem}
\begin{proof}
By the convergence of eigenvalues ensured by Lemma
\ref{l:conv_of_eigenvalues}, ${\mathcal D}$ is an open set. Therefore, by
Proposition \ref{p:open-dense_direct_sum}, ${\mathcal M} \cap
{\mathcal D}$ is also an open set. The geometric statement is straightforward. Furthermore, since ${\mathcal M}$ is dense
in $M(n,\bbR)$, then ${\mathcal M} \cap {\mathcal D}$ must also be
dense in the relative topology of the open set ${\mathcal D}$.

The converse is an immediate consequence of Proposition
\ref{p:open-dense_direct_sum}. $\Box$
\end{proof}

\begin{remark}\label{r:heuristic_proof_identifiability}
A simple heuristic argument may shed light on a weaker version of the claim of Theorem \ref{t:minimal_is_topol_general}, i.e., with identifiability (uniqueness) of the matrix exponent $H$ in place of the minimality of the symmetry group.

Consider a matrix parameter $S+L = M = W^{-1}DW$ such that $S \in {\mathcal S}_{\neq}$ and $L \in so(n)$. For $x \neq 1$, but close to 1, the Baker-Campbell-Hausdorff formula gives
$$
\textnormal{Log}(x^{-M}x^{-M^*}) = -\log(x) (M+M^*) + O(\log^2(x)).
$$
Since
$$
\frac{O(\log^2(x))}{\log(x)}\rightarrow {\mathbf 0}, \quad x \rightarrow 1,
$$
then the term $O(\log^2(x))$ is only a slight perturbation of the matrix $-\log(x) (M+M^*)$. Assume for simplicity that $O(\log^2(x))\equiv{\mathbf 0}$. Then $\Pi_x = e^{\textnormal{Log}(x^{-M}x^{-M^*})}= e^{-\log(x)(M+M^*)}$ is symmetric with pairwise different eigenvalues. We thus obtain
$$
G(\Pi_x) = P_x (\pm 1, \pm 1, \hdots, \pm1 )P^{*}_{x},
$$
where $P_x \in O(n)$ is a matrix of eigenvectors of $\Pi_x$. In particular, $G(\Pi_x)$ is finite and thus the matrix exponent of the associated OFBM is identifiable.

Moreover, the set ${\mathcal S}_{\neq}$ is open and dense in $S(n,\bbR)$, thus implying that ${\mathcal S}_{\neq} \oplus so(n)$ is open and dense in the subset of $M(n,\bbR)$ whose eigenvalues have real parts between $-1/2$ and 1/2.
\end{remark}

%
%\begin{corollary}
%For $n \geq 2$, minimal symmetry holds except on a meager set of matrix parameters $M$.
%\end{corollary}
%\begin{proof}
%Choose $M = S+L$, where $S \in {\mathcal S}_{\neq}$ and $L \in ({\mathcal L}_{\textnormal{invar}}(o_1,\hdots,o_n))^c$, where $o_1,\hdots,o_n$ are a basis of orthonormal eigenvectors of $S$. Since ${\mathcal S}_{\neq}$ and $({\mathcal L}_{\textnormal{invar}}(o_1,\hdots,o_n))^c$ are open and dense in their relative topologies, then the graph of the function $S \mapsto ({\mathcal L}_{\textnormal{invar}}(o_1,\hdots,o_n))^c$ is also open and dense in the topology of $M(n,\bbR)$ (by the same argument as in \texttt{ofbm2}). $\Box$
%\end{proof}

\section{Classification in dimensions $n=2$ and $n=3$}
\label{s:2-3}

Theorems \ref{t:symmetry-group-OFBM} and
\ref{t:symmetry-group-OFBM-2} describe the general structure of
symmetry groups of OFBMs. The cases of maximal and minimal symmetry groups were studied in
Section \ref{s:min-max}. In this section, we are interested in identifying all the possible
``intermediate" symmetry groups. We shall describe their structure in
dimensions $n=2$ and $n=3$, and make some comments about higher dimensions.

\subsection{Dimension $n=2$}
\label{s:2}

When $n=2$, the contribution of the term $G(\Pi_I)$ in Theorems
\ref{t:symmetry-group-OFBM} and \ref{t:symmetry-group-OFBM-2} can be
easily described, as the next two results show.

\begin{lemma}\label{l:2-Pi-I}
If $\Pi_I \neq 0$, then $G(\Pi_I) = SO(2)$.
\end{lemma}

\begin{proof}
Since $\Pi_I$ is skew-symmetric, we have
$$
\Pi_I = \lambda \left(\begin{array}{cc} 0 & 1 \\
-1 & 0
\end{array}
\right),\quad \lambda \neq 0.
$$
Thus, $\Pi_I/\lambda$ is a rotation matrix, and thus
$G(\Pi_{I})=SO(2)$. \ \ $\Box$
\end{proof}

\medskip

Theorems \ref{t:symmetry-group-OFBM} and
\ref{t:symmetry-group-OFBM-2} can now be reformulated as follows.

\begin{corollary}\label{c:2-structure}
For $n=2$, under the assumptions and notation of Theorems
\ref{t:symmetry-group-OFBM} and \ref{t:symmetry-group-OFBM-2}, we have
\begin{eqnarray}
% \nonumber to remove numbering (before each equation)
  G_H &=& W \left\{\begin{array}{cl}
    \cap_{x>0} G(\Pi_x) \cap SO(2), & \mbox{\rm if}\ \Im(AA^*) \neq 0\\
    \cap_{x>0} G(\Pi_x), & \mbox{\rm if}\ \Im(AA^*) = 0
  \end{array} \right\} W^{-1}
  \label{e:2-structure} \\
   &=& W \left\{\begin{array}{cl}
    \cap_{m\geq 1} G(\Pi^{(m)}) \cap SO(2), & \mbox{\rm if}\ \Im(AA^*) \neq 0\\
    \cap_{m\geq 1} G(\Pi^{(m)}), & \mbox{\rm if}\ \Im(AA^*) = 0
  \end{array} \right\} W^{-1}.
   \label{e:2-structure-2}
\end{eqnarray}

\end{corollary}

Next, we study the possible structures of the groups $G(\Pi)$ when
$\Pi$ is symmetric (and hence potentially positive definite, as the matrix $\Pi_x$ in
(\ref{e:2-structure})). Let $\pi_1,\pi_2$ be the two real
eigenvalues of $\Pi$. Two cases need to be considered:
\begin{equation}\label{e:2-cases}
    \begin{array}{cc}
    \mbox{Case 2.1:} & \pi_1 = \pi_2,\\
    \mbox{Case 2.2:} & \pi_1 \neq \pi_2.
    \end{array}
\end{equation}
In Case 2.1, $\Pi = \pi_1 I$ and hence
\begin{equation}\label{e:2-case2.1}
    G(\Pi) = O(2).
\end{equation}
In Case 2.2, we can write
$$
\Pi = S \left( \begin{array}{cc}
\pi_1 & 0 \\
0 & \pi_2
\end{array}
\right) S^* = (p_1 , p_2) \left( \begin{array}{cc}
\pi_1 & 0 \\
0 & \pi_2
\end{array}
\right) (p_1 , p_2)^{*},
$$
where the columns of the orthogonal matrix $S = (p_1\ p_2)$ consist of the orthonormal
eigenvectors $p_1,p_2$ of $\Pi$. By Theorem \ref{t:comm_general}, $B\in O(2)$ commutes with such $\Pi$ if and only if $B =
SGS^*$ where $G$ is a diagonal matrix such that $G^2 = I$ ($G^2 =I$ is a
consequence of the fact that $B\in O(2)$), or
$$
B = S \left( \begin{array}{cc}
\pm 1 & 0 \\
0 & \pm 1
\end{array}
\right) S^* = (p_1, p_2) \left( \begin{array}{cc}
\pm 1 & 0 \\
0 & \pm 1
\end{array}
\right) (p_1 , p_2)^{*}.
$$
We thus have
\begin{eqnarray}
% \nonumber to remove numbering (before each equation)
 G(\Pi)  &=& \Big\{ I,-I, S \left( \begin{array}{cc}
1 & 0 \\
0 & - 1
\end{array}
\right) S^*, S \left( \begin{array}{cc}
- 1 & 0 \\
0 & 1
\end{array}
\right) S^*  \Big \}  \nonumber \\
   &=& \{I,-I,\mbox{Ref}(p_1),\mbox{Ref}(p_2) \},
\label{e:2-case2.2}
\end{eqnarray}
where $\mbox{Ref}(p)$ indicates a reflection around the axis
spanned by a vector $p$. The expressions (\ref{e:2-case2.1}) and
(\ref{e:2-case2.2}) provide the only possible structures for
$G(\Pi)$. Together with Corollary \ref{c:2-structure}, this leads
to the following result.

\begin{theorem}\label{t:2-structure}
Consider an OFBM given by the spectral representation
(\ref{e:spectral-repres-OFBM}), and suppose that the matrix $A$
satisfies the assumption (\ref{e:full-rank}). Then, its symmetry
group $G_H$ is conjugate to one of the following:
\begin{itemize}
    \item [$(2.a)$] minimal type: $\{I,-I\}$;
    \item [$(2.b)$] trivial type: $\{I,-I,\mbox{\rm Ref}(p_1),\mbox{\rm Ref}(p_2)\}$ for a pair of
    orthogonal $p_1,p_2$;
    \item [$(2.c)$] rotational type: $SO(2)$;
    \item [$(2.d)$] maximal type: $O(2)$.
\end{itemize}
\end{theorem}

All the types of subgroups described in Theorem \ref{t:2-structure} are non-empty, as we show next.
Since OFBMs of maximal and minimal types were studied in general dimension $n$ in Section
\ref{s:min-max}, we now provide examples of OFBMs of only the two remaining types for dimension $n=2$.

\begin{example}\label{ex:2-rotational} (Rotational type)
Consider an OFBM with parameters
\begin{equation}\label{e:2-rotational}
    D = dI,\quad \sqrt{2}A_1 \in SO(2) \backslash\{I,-I\}, \quad \sqrt{2}A_2 \in O(2)\backslash SO(2),
\end{equation}
where $d$ is real. Then, we have $\Pi_x = x^{-2d}I$
and $G(\Pi_x)=O(2)$. Since $\Im(AA^*) \neq 0$, Corollary
\ref{c:2-structure} yields that $G_H = SO(2)$.
\end{example}

\begin{example}\label{ex:2-trivial} (Trivial type)
Consider an OFBM with parameters
\begin{equation}\label{e:2-trivial}
    D = \left(\begin{array}{cc} d_1&0\\ 0&d_2\end{array}\right),\quad
    A = \left(\begin{array}{cc} a_1&0\\ 0&a_2\end{array}\right),
\end{equation}
where $d_1\neq d_2$ are real. Then,
$$
AA^* = \left(\begin{array}{cc} |a_1|^2 &0\\
0&|a_2|^2\end{array}\right) = \Re(AA^*),\quad \Im(AA^*) = 0,
$$
and
$$
\Pi_x = \left(\begin{array}{cc} x^{-2d_1}&0\\
0&x^{-2d_2}\end{array}\right),
$$
implying that, for $x\neq 1$,
$$
G(\Pi_x) = \left\{I,-I,\left(\begin{array}{cc} 1&0\\
0&-1\end{array}\right),\left(\begin{array}{cc} -1&0\\
0&1\end{array}\right) \right\}.
$$
Corollary \ref{c:2-structure} then yields
\begin{equation}\label{e:2-trivial-G_H}
    G_H =\left\{I,-I,\left(\begin{array}{cc} 1&0\\
0&-1\end{array}\right),\left(\begin{array}{cc} -1&0\\
0&1\end{array}\right) \right\}.
\end{equation}
\end{example}

Only OFBMs of rotational and maximal types have multiple
exponents. Moreover, in view of (\ref{e:exponents}), in both cases
we have
\begin{equation}\label{e:2-nonunique}
    {\cal E}(B_H) = H + Wso(2)W^{-1},
\end{equation}
where $H$ is any exponent of the OFBM $B_H$. This relation can be
further refined, as the following proposition shows. For this purpose, we need to consider
a so-called commuting exponent $H_0 \in {\mathcal E}(B_H)$, i.e., an exponent $H_0$ such that
\begin{equation}\label{e:H0_def}
H_0 C = C H_0
\end{equation} for all $C\in G_H$. The existence of this useful exponent is ensured by Lemma 2 of Maejima
\cite{maejima:1998}.

\begin{proposition}\label{p:2-exponents}
Consider an OFBM given by the spectral representation
(\ref{e:spectral-repres-OFBM}), and suppose that the matrix $A$
satisfies the assumption (\ref{e:full-rank}). If ${\mathcal E}(B_H)$ is not unique, then
the commuting exponents are of the form
\begin{equation} \label{e:2-H0_nonunique_E}
H_0 = W U_2 \textnormal{diag}(h,\overline{h}) U^{*}_2 W^{-1},
\end{equation}
where
\begin{equation}\label{e:U2}
U_2 = \frac{\sqrt{2}}{2}\left(\begin{array}{cc}
1 & 1 \\
i & -i
\end{array}\right)
\end{equation}
and $h \in \bbC$. In
particular,
\begin{equation} \label{e:2-nonunique_E}
{\mathcal E}(B_H) = W(U_2 \textnormal{diag}(h,\overline{h})U^{*}_2
+ so(2))W^{-1},
\end{equation}
$H = \Re(h)I \in {\mathcal E}(B_H)$ and $W^{-1}H W$ is
normal for any $H \in {\mathcal E}(B_H)$.
\end{proposition}

\begin{proof}
If ${\mathcal E}(B_H)$ is not unique, then by Theorem
\ref{t:2-structure}, $H_0$ commutes with $W SO(2)W^{-1}$. In
particular, $H_0$ commutes with $W O W^{-1}$ for $O \in SO(2)
\backslash\{I,-I\}$. Since such $O$ is diagonalizable with two
complex conjugate eigenvalues, the eigenvectors of $W O W^{-1}$
are also eigenvectors of $H_0$. Thus, $H_0$ can be written as $W
U_2 \textnormal{diag}(h_1,h_2) U^{*}_2W^{-1}$. Therefore, since
$h_1$, $h_2$ are also the eigenvalues of $U_2
\textnormal{diag}(h_1,h_2) U^{*}_2$, which must only have real
entries, a simple calculation shows that $h_1$ = $\overline{h}_2$,
and thus (\ref{e:2-H0_nonunique_E}) holds. This also yields
(\ref{e:2-nonunique_E}).

For $H \in {\mathcal E}(B_H)$, (\ref{e:2-nonunique_E}) implies
that $W^{-1}HW$ is normal. In particular, we may choose the
exponent $H = H_0 + W
L_{-\Im(h)}W^{-1}= \Re(h)I$, where $L_s$ is
defined in (\ref{e:L_s}). $\Box$
\end{proof}

\begin{remark}
In the case of OFBMs of rotational type, which have multiple exponents, every exponent is a commuting exponent
(compare with Meerschaert and Veeh \cite{meerschaert:veeh:1993}, p.\ 721, for the case of operator stable measures).

\end{remark}

\begin{remark}
For general proper Gaussian processes, one can define symmetry sets (groups) in
the same way as for o.s.s.\ processes, and, in particular, show
that they are also compact subgroups of $GL(n,\bbR)$. By applying
the argument of the proof of Theorem 4.5.3 in Didier
\cite{didier:2007}, which is based on general commutativity
results, instead of spectral filters, one can show that the
classification provided by Theorem \ref{t:2-structure} actually
holds for the wide class of proper bivariate Gaussian processes.
\end{remark}

\subsection{Dimension $n=3$}
\label{s:3}

We will make use of the partition of $O(3)$ into the
following subsets:
$$
SO(3) = \{I\} \cup \textnormal{Rot}_\theta \cup
\textnormal{Rot}_\pi, \quad
O(3) \backslash SO(3) = \{-I\} \cup \textnormal{Ref}_\theta
\cup \textnormal{Ref}_0,
$$
where for a vector $p$,
$$
\textnormal{Rot}_\theta  := \bigcup_{p \in S^{n-1}} \textnormal{Rot}_\theta(p), \quad \textnormal{Rot}_\pi  := \bigcup_{p \in S^{n-1}} \textnormal{Rot}_\pi(p),
$$
$$
\textnormal{Ref}_\theta  := \bigcup_{p \in S^{n-1}} \textnormal{Ref}_\theta(p), \quad \textnormal{Ref}_0  := \bigcup_{p \in S^{n-1}} \textnormal{Ref}_0(p),
$$
and
\begin{eqnarray}
% \nonumber to remove numbering (before each equation)
  \mbox{Rot}_\theta(p) &=& \{\mbox{rotations about the axis $\textnormal{span}_{\bbR}(p)$ by an angle not equal to $\pi$}\}, \nonumber  \\
  \mbox{Rot}_\pi(p) &=& \{\mbox{rotation about the axis $\textnormal{span}_{\bbR}(p)$ by an angle equal to $\pi$}\}, \nonumber  \\
  \mbox{Ref}_\theta(p) &=&  \{\mbox{rotations about the axis $\textnormal{span}_{\bbR}(p)$ by an angle not equal to $\pi$, combined with the }  \nonumber  \\
  & & \mbox{reflection in the plane through the origin which is perpendicular to the axis}\}, \nonumber  \\
  \mbox{Ref}_0(p) &=& \{\mbox{reflection in a plane through the origin, where the plane is perpendicular to $p$}\}. \nonumber  \\
 % & & \mbox{reflection in the plane through the origin which is perpendicular to the axis}\}. \nonumber
\end{eqnarray}

From a matrix perspective, for some $p \in S^{n-1}$,
$$
\textnormal{Rot}_\theta(p) \cong \left(\begin{array}{ccc} \cos\theta &
-\sin \theta & 0\\
\sin \theta & \cos \theta & 0\\
0 & 0 & 1
\end{array}\right), \hspace{2mm}\theta \in (0,2\pi)\backslash \{\pi\}, \quad \textnormal{Rot}_\pi(p) \cong \left(\begin{array}{ccc} -1 &
0 & 0\\
0 & -1 & 0\\
0 & 0 & 1
\end{array}\right),
$$
$$
\textnormal{Ref}_\theta(p) \cong \left(\begin{array}{ccc} \cos\theta &
-\sin \theta & 0\\
\sin \theta & \cos \theta & 0\\
0 & 0 & -1
\end{array}\right), \hspace{2mm} \theta \in (0,2\pi)\backslash \{\pi\},
\quad \textnormal{Ref}_0 (p) \cong \left(\begin{array}{ccc} 1 &
0 & 0\\
0 & 1 & 0\\
0 & 0 & -1
\end{array}\right),
$$
where $\cong$ indicates conjugacies by orthogonal matrices.

\begin{remark}
The subscript $\theta$ in $\textnormal{Rot}_\theta$ or
$\textnormal{Ref}_\theta$ only indicates that the angle in
question is not 0 or $\pi$. Here, $\theta$ does \textit{not} refer
to a specific angle. Indeed, even in the case of a fixed $p$,
$\textnormal{Rot}_\theta(p)$ and $\textnormal{Ref}_\theta(p)$ are
classes of matrices. Also, in the expression $\textnormal{Ref}_0$
we use the subscript 0 to indicate that there is no rotation
before reflection through the plane in question.
\end{remark}

We first describe the possible structures of $G(\Pi)$ for
symmetric matrices $\Pi$ (such as the matrices $\Pi_x$, $x > 0$, in
(\ref{e:symmetry-group-OFBM-main})). Let $\pi_1,\pi_2,\pi_3$ be
the three real eigenvalues of $\Pi$. Three cases need to be
considered, namely,
\begin{equation}\label{e:3-cases}
    \begin{array}{cl}
    \mbox{Case 3.1:} & \pi_1 = \pi_2 = \pi_3 ,\\
    \mbox{Case 3.2:} & \pi_1 = \pi_2 \neq \pi_3, \\
    \mbox{Case 3.3:} & \pi_i \neq \pi_j, \hspace{2mm} i \neq j, \hspace{2mm} i,j=1,2,3.
    \end{array}
\end{equation}

The next proposition gives the form of $G(\Pi)$ in all the above cases.

\begin{proposition}\label{p:3-G(Pi)}
Let $\Pi \in {\mathcal S}(3,\bbR)$. Denote its eigenvectors by $p_i$, $i=1,2,3$, where $S = (p_1 \ p_2 \ p_3) \in O(3)$. Then,
\begin{enumerate}
\item [(i)] in Case 3.1 in (\ref{e:3-cases}),
\begin{equation}\label{e:3-case3.1}
    G(\Pi) = O(3);
\end{equation}
\item [(ii)] in Case 3.2 in (\ref{e:3-cases}),
$$
G(\Pi) = \{I,-I\} \cup (\textnormal{Rot}_\theta(p_3) \cup
\textnormal{Ref}_\theta(p_3))
     \cup (\textnormal{Rot}_\pi(p_3) \cup \textnormal{Ref}_0(p_3))
$$
\begin{equation}\label{e:3-case3.2}
     \cup
     \bigcup_{q\in \mbox{\scriptsize span}_{\bbR}\{p_1,p_2\}} (\textnormal{Rot}_\pi(q) \cup
     \textnormal{Ref}_0(q));
\end{equation}
\item [(iii)] in Case 3.3 in (\ref{e:3-cases}),
\begin{equation}\label{e:3-case3.3}
G(\Pi) = \{I,-I,\textnormal{Ref}_0(p_1),\textnormal{Ref}_0(p_2),\textnormal{Ref}_0(p_3),
\textnormal{Rot}_\pi(p_1),\textnormal{Rot}_\pi(p_2),\textnormal{Rot}_\pi(p_3)\}.
\end{equation}
\end{enumerate}
\end{proposition}
\begin{proof}
$(i)$ is immediate, so we turn to $(ii)$. In this case, we can write $\Pi = S\textnormal{diag}(\pi_1,\pi_1,\pi_3)S^*$. By Theorem
\ref{t:comm_general}, $B$ commutes with such $\Pi$ if and only if
\begin{equation}\label{e:matrix_B}
B = S \left( \begin{array}{ccc}
c_{11} & c_{12} & 0\\
c_{21} & c_{22} & 0 \\
0 & 0 & d
\end{array}
\right) S^*,
\end{equation}
where $C = (c_{ij})_{i,j=1,2}$ and $d$ are arbitrary. If we are
only interested in orthogonal matrices, this gives $C \in O(2)$, and $d=\pm 1$, which corresponds
to the subgroup (\ref{e:3-case3.2}).
Indeed, the matrices $\textnormal{Rot}_{\theta}(p_3)$ and
$\textnormal{Rot}_{\pi}(p_3)$ in (\ref{e:3-case3.2}) account for
rotations $C \in SO(2)$ and $d=1$ in (\ref{e:matrix_B}),
$\textnormal{Ref}_{\theta}(p_3)$ and $\textnormal{Ref}_{0}(p_3)$
in (\ref{e:3-case3.2}) account for rotations $C \in SO(2)$ and $d=-1$ in
(\ref{e:matrix_B}), and $\textnormal{Rot}_{\pi}(q)$ and
$\textnormal{Ref}_{0}(q)$, $q \in
\textnormal{span}_{\bbR}\{p_1,p_2\}$, account for reflections $C \in O(2)\backslash SO(2)$
and $d = \pm 1$ in (\ref{e:matrix_B}).

Regarding $(iii)$, we can write $\Pi = S\textnormal{diag}(\pi_1,\pi_2,\pi_3)S^*$.
By Theorem \ref{t:comm_general}, $B\in O(3)$ commutes with such $\Pi$ if and only if $B = S \textnormal{diag}(\pm 1,\pm 1,\pm 1)S^*$.
We thus have
\begin{eqnarray}
G(\Pi) & = &  \{I,-I, S \textnormal{diag}(-1,1,1)S^*,S
\textnormal{diag}(1,-1,1)S^*,S \textnormal{diag}(1,1,-1)S^*, \nonumber \\
& & S \textnormal{diag}(1,-1,-1)S^*,S
\textnormal{diag}(-1,1,-1)S^*,S \textnormal{diag}(-1,-1,1)S^*\}, \nonumber
\end{eqnarray}
as stated. $\Box$\\
\end{proof}

The expressions (\ref{e:3-case3.1}), (\ref{e:3-case3.2}) and
(\ref{e:3-case3.3}) describe the only possible structures for
$G(\Pi)$, all of them, as shown below, being symmetry groups
of some non-empty subclass of OFBMs. However, new symmetry groups
may arise when one considers intersections of $G(\Pi_{x})$ for different
values of $x$, and also with $G(\Pi_{I})$. In order to provide
a full description of symmetry groups of OFBMs in dimension $n=3$, we first
consider the case of time reversible OFBMs, before turning to the general case.
As shown in Didier and Pipiras
\cite{didier:pipiras:2009}, time reversibility corresponds to the assumption that
\begin{equation}\label{e:time-reversible}
    \Im(AA^*) = 0.
\end{equation}
Under (\ref{e:time-reversible}), the
presence of $G(\Pi_I)$ in (\ref{e:symmetry-group-OFBM-main}) and
(\ref{e:symmetry-group-OFBM-main-2})
can be ignored.

\begin{theorem}\label{t:3-structure}
Consider an OFBM given by the spectral representation
(\ref{e:spectral-repres-OFBM}), and suppose that the matrix $A$
satisfies the assumptions (\ref{e:full-rank}) and
(\ref{e:time-reversible}). Then, its symmetry group $G_H$ is
conjugate by a positive definite matrix $W$ to one of the
following:
\begin{itemize}
    \item [$(3.a)$] minimal type: $\{I,-I\}$;
    \item [$(3.b)$]  for some vector $p$,
$$
\{I,-I,\mbox{\rm Ref}_0(p),\mbox{\rm Rot}_\pi(p)\};
$$
    \item [$(3.c)$] for some orthogonal $p_1,p_2,p_3$,
$$
\{I,-I,\mbox{\rm Ref}_0(p_1),\mbox{\rm Ref}_0(p_2),\mbox{\rm
Ref}_0(p_3), \mbox{\rm Rot}_\pi(p_1),\mbox{\rm
Rot}_\pi(p_2),\mbox{\rm Rot}_\pi(p_3)\};
$$
    \item [$(3.d)$] for some orthogonal $p_1,p_2,p_3$,
$$
\{I,-I\} \cup (\mbox{\rm Rot}_\theta(p_3) \cup \mbox{\rm
Ref}_\theta(p_3))
     \cup (\mbox{\rm Rot}_\pi(p_3) \cup \mbox{\rm Ref}_0(p_3))     \cup
     \bigcup_{q\in \mbox{\scriptsize span}_{\bbR}\{p_1,p_2\}} (\mbox{\rm Rot}_\pi(q) \cup
     \mbox{\rm Ref}_0(q)) ;
$$
    \item [$(3.e)$] maximal type: $O(3)$.
\end{itemize}
\end{theorem}
\begin{proof}
Recall that, under (\ref{e:time-reversible}), the symmetry group
$G_H$ is conjugate to $\cap_{x > 0}G(\Pi_x)$. By Proposition \ref{p:3-G(Pi)}, $G(\Pi_x)$ can only be of the forms
(\ref{e:3-case3.1}), (\ref{e:3-case3.2}) and (\ref{e:3-case3.3}).
The proof is now split into the following cases.\\

Case 1: for some $x > 0$, $G(\Pi_x)$ has the form
(\ref{e:3-case3.3}). Since the intersection with some other
$G(\Pi_{x'})$ can only reduce the group, $\cap_{x > 0}G(\Pi_x)$
can be only of types $(3.a)$, $(3.b)$ or $(3.c)$, where $(3.b)$ is
a consequence of intersecting (\ref{e:3-case3.3}) with some
appropriate (\ref{e:3-case3.3}).\\

Case 2: all $G(\Pi_x)$ are the same and have the form
(\ref{e:3-case3.2}). This gives type $(3.d)$.\\

Case 3: there are $x_1 \neq x_2$ and two different $G(\Pi_{x_1})$
and $G(\Pi_{x_2})$ such that both have the form
(\ref{e:3-case3.2}). Let $r_1, r_2, r_3$ and
$p_1, p_2, p_3$  be the corresponding orthonormal
vectors in (\ref{e:3-case3.2}). For $G(\Pi_{x_1})$ and
$G(\Pi_{x_2})$ to be different, the corresponding axis
$\textnormal{span}_{\bbR}\{r_3\}$ and
$\textnormal{span}_{\bbR}\{p_3\}$ have to be different. We break up the proof into subcases. In all of them,
based on necessary conditions for commutativity, we will construct a set ${\mathcal C}$ of ``candidate" matrices such that ${\mathcal C} \supseteq G(\Pi_{x_1})\cap G(\Pi_{x_2})$. As we will see, in all subcases, ${\mathcal C}$ has one of the forms $(3.a)$, $(3.b)$ or $(3.c)$. Since all the possible symmetry groups which are subgroups of $(3.c)$ (the most encompassing of the three) are of the forms $(3.a)$, $(3.b)$ or $(3.c)$, then $\cap_{x > 0} G(\Pi_{x})$ must be of one of the latter forms, which completes the proof.

   Thus, the subcases and the respective sets ${\mathcal C}$ are as follows.

\begin{itemize}
\item[$(3.i)$] $p_3 \in \textnormal{span}_{\bbR}\{r_1,r_3\} \backslash (\textnormal{span}_{\bbR}\{r_1\} \cup \textnormal{span}_{\bbR}\{r_3\})$ (or with $r_2$ in place of $r_1$): then, by Theorem \ref{t:comm_general} (c.f.\ the proof of Proposition \ref{p:3-G(Pi)}, case $(ii)$), both $\textnormal{span}_{\bbR}\{r_3\}$ and $\textnormal{span}_{\bbR}\{p_3\}$ are real invariant subspaces of any $M \in G(\Pi_{x_1}) \cap G(\Pi_{x_2})$. Consequently, since $r_3$ and $p_3$ are not orthogonal, and two real eigenvectors of orthogonal matrices associated with different eigenvalues 1 and $-1$ must be orthogonal, for $M$ we must have that $\textnormal{span}_{\bbR}\{r_3,p_3\}=\textnormal{span}_{\bbR}\{r_1,r_3\}$ is a two-dimensional real (proper or not) subspace of the eigenspace associated with either 1 or $-1$. Therefore, the remaining eigenvalue of $M$ must be real, and by Lemma \ref{l:third_eigenvec_of_O} there exists an associated eigenvector which is orthogonal to $\textnormal{span}_{\bbR}\{r_1,r_3\}$. In particular, $r_2$ is an eigenvector. Therefore, we obtain the set of candidate matrices
$$
{\mathcal C} = \{\pm I, \pm (r_1,r_2,r_3)\textnormal{diag}(1,-1,1)(r_1,r_2,r_3)^{*}\}.
$$

\item[$(3.ii)$] $p_3 \in \textnormal{span}_{\bbR}\{r_1,r_2\}$: then for $M \in G(\Pi_{x_1}) \cap G(\Pi_{x_2})$, by the same argument as in $(3.i)$, both $p_3 \perp r_3$ are (real) eigenvectors of $M$ associated with real eigenvalues. Thus, the third eigenvalue of $M$ is also real and, by Lemma \ref{l:third_eigenvec_of_O}, $(\textnormal{span}_{\bbR}\{p_3,r_3\})^{\perp} \subseteq \textnormal{span}_{\bbR}\{r_1,r_2\}$ contains a unit norm real eigenvector of $M$, which we can denote by $q$. Thus, we can set
    $$
    {\mathcal C} = \{M: M = (q,p_3,r_3)\textnormal{diag}(\pm 1,\pm 1,\pm 1)(q,p_3,r_3)^{*}\}.
    $$

   \item[$(3.iii)$] $p_3 \notin \cup_{i \neq j}\textnormal{span}_{\bbR}\{r_i,r_j\}$: then, by the same argument as in $(3.i)$,
$\textnormal{span}_{\bbR}\{r_3\}$ and $\textnormal{span}_{\bbR}\{p_3\}$ are both real invariant subspaces for the solutions.
However, since they are not orthogonal, then, for each $M \in G(\Pi_{x_1})\cap G(\Pi_{x_2})$, we have that $\textnormal{span}_{\bbR}\{p_3,r_3\}$ is
a two-dimensional real (proper or not) subspace of the eigenspace associated with either 1 or $-1$. Thus, by Lemma \ref{l:third_eigenvec_of_O},
there exists $q \in \textnormal{span}_{\bbR}\{p_3,r_3\}^{\perp}$, $\norm{q}=1$, such that we can write
       $$
       {\mathcal C} = \{ M:  M = (q,p_3,r_3)\textnormal{diag}(\pm 1, \pm I_{2})(q,p_3,r_3)^{-1} \}.
       $$

       Note that the matrices in ${\mathcal C}$ are, indeed, orthogonal. In fact, set $v_1 = q$ and $v_3 = r_3$. Now let $\pi_{r_1,r_2}$ be the operator for the projection on $\textnormal{span}_{\bbR}\{r_1,r_2\}$, and set $v_2 = \frac{\pi_{r_1,r_2}(p_3)}{ \norm{\pi_{r_1,r_2}(p_3)} }$. Then $v_2 \neq {\mathbf 0}$, because $p_3$ and $r_3$ are not collinear. Moreover, we have
       $$
       \langle v_2, q \rangle= \frac{1}{\norm{\pi_{r_1,r_2}(p_3)}}\langle\pi_{r_1,r_2}(p_3), q \rangle = \frac{1}{\norm{\pi_{r_1,r_2}(p_3)}}\langle p_3, \pi_{r_1,r_2}(q) \rangle= \frac{1}{\norm{\pi_{r_1,r_2}(p_3)}} \langle p_3, q \rangle = 0,
       $$
   where we used the self-adjointness of $\pi_{r_1,r_2}(\cdot)$ and the fact that $q \in \textnormal{span}_{\bbR}\{r_1,r_2\}$. Therefore, we have that $v_2 \in \textnormal{span}_{\bbR}\{p_3, v_3\}$. Moreover, since $v_2 \perp r_3$, then $v_1, v_2, v_3$ are linearly independent and thus form an orthonormal basis. Then we can write
   $$
       {\mathcal C} = \Big\{ M:  M = \Big(q,\frac{\pi_{r_1,r_2}(p_3)}{ \norm{\pi_{r_1,r_2}(p_3)} },r_3 \Big)\textnormal{diag}(\pm 1, \pm I_{2})\Big(q,\frac{\pi_{r_1,r_2}(p_3)}{ \norm{\pi_{r_1,r_2}(p_3)} },r_3 \Big)^{*} \Big\}.
       $$

\end{itemize}

Case 4: For all $x > 0$, $G(\Pi_x)$ has the form
(\ref{e:3-case3.1}). This gives type $(3.e)$. $\Box$\\
\end{proof}

We now provide examples of OFBMs of the types
$(3.b)$, $(3.c)$ and $(3.d)$, thereby showing that all the types described in Theorem \ref{t:3-structure} are non-empty.

\begin{example}\label{ex:3-3.b} (Type $(3.b)$)
Consider the OFBM with spectral representation parameters $A :=I$ and
$$
D = \left(\begin{array}{ccc}
d & 0 & 0\\
1 & d & 0 \\
0 &0 & d
\end{array}\right).
$$
By Theorem \ref{t:symmetry-group-OFBM}, we may assume that the
positive definite conjugacy associated with $G(B_H)$ is $W = I$.
Observe that
$$
\Pi_x = x^{-D}x^{-D^*} = x^{-2d}\left(\begin{array}{ccc}
1 & -\log(x) & 0\\
-\log(x) & \log^2(x)+1 & 0\\
0 &0 & 1
\end{array}\right).
$$
Due to the block-diagonal shape of $x^{-D}x^{-D^*}$, it suffices
to focus on its $2 \times 2$ upper left block. Consider $x =
e^{-1}, e$. The associated $2 \times 2$ blocks, i.e.,
$$
\left(\begin{array}{cc}
1 & 1 \\
1 & 2
\end{array}\right), \quad
\left(\begin{array}{cc}
1 & -1\\
-1 & 2\\
\end{array}\right),
$$
have pairwise different eigenvalues. Moreover, they do not share eigenvectors, since
otherwise they would commute. As a consequence, by Proposition \ref{p:3-G(Pi)},
 they are of the form $(3.b)$ with $p= (0,0,1)'$.
\end{example}

\begin{example}\label{ex:3-3.c} (Type $(3.c)$)
Consider the OFBM with spectral representation parameters
\begin{equation}\label{e:3-3.c}
    D = \mbox{diag}(d_1,d_2,d_3),\ A = \mbox{diag}(a_1,a_2,a_3),
\end{equation}
where $d_i \neq d_j$, $i \neq j$. Then, $AA^* =
\mbox{diag}(|a_1|^2,|a_2|^2,|a_3|^2) = \Re(AA^*)$, $\Im(AA^*) = 0$
and $\Pi_x = \mbox{diag}(x^{-2d_1},x^{-2d_2},x^{-2d_3})$. This
yields
\begin{equation}\label{e:3-3.c-G_H}
    G_H = \{I,-I,\mbox{\rm Ref}_0(e_1),\mbox{\rm Ref}_0(e_2),\mbox{\rm
Ref}_0(e_3), \mbox{\rm Rot}_\pi(e_1),\mbox{\rm
Rot}_\pi(e_2),\mbox{\rm Rot}_\pi(e_3)\}
\end{equation}
for the Euclidean vectors $e_i$, $i=1,2,3$.
\end{example}

\begin{example}\label{ex:3-3.d} (Type $(3.d)$)
Consider the OFBM with spectral representation parameters
\begin{equation}\label{e:3-3.d}
    D = \mbox{diag}(d_1,d_1,d_3),\ A = \mbox{diag}(a_1,a_2,a_3),
\end{equation}
where $d_1\neq d_3$. Then, $AA^* =
\mbox{diag}(|a_1|^2,|a_2|^2,|a_3|^2) = \Re(AA^*)$, $\Im(AA^*) = 0$
and $\Pi_x = \mbox{diag}(x^{-2d_1},x^{-2d_1},x^{-2d_3})$. This
yields
\begin{equation}\label{e:3-3.d-G_H}
    G_H = \{I,-I\} \cup (\mbox{\rm Rot}_\theta(e_3) \cup \mbox{\rm
Ref}_\theta(e_3))
     \cup (\mbox{\rm Rot}_\pi(e_3) \cup \mbox{\rm Ref}_0(e_3))     \cup
     \bigcup_{q\in \mbox{\scriptsize span}\{e_1,e_2 \}} (\mbox{\rm Rot}_\pi(q) \cup
     \mbox{\rm Ref}_0(q)).
\end{equation}
\end{example}

We now extend Theorem \ref{t:3-structure} to the general case of
OFBMs which are not necessarily time reversible, i.e., we drop the assumption (\ref{e:time-reversible}). From the perspective of the structural result provided by Theorem \ref{t:symmetry-group-OFBM}, the lack of time reversibility manifests itself as an additional constraint which may reduce the symmetry group, and even generate a new type, as seen in the next theorem.

\begin{theorem}\label{t:3-structure_nontimerevers}
Consider an OFBM given by the spectral representation
(\ref{e:spectral-repres-OFBM}), and suppose that the matrix $A$
satisfies the assumption (\ref{e:full-rank}). Then, its symmetry
group $G_H$ is conjugate by a positive definite matrix $W$ to the
ones described in Theorem \ref{t:3-structure}, plus the following:
\begin{itemize}
\item [$(3.f)$] for some vector $p$,
$$
\{I,-I,\mbox{\rm Ref}_0(p),\mbox{\rm Rot}_\pi(p),\mbox{\rm
Ref}_\theta(p),\mbox{\rm Rot}_\theta(p)\}.
$$
\end{itemize}
\end{theorem}
\begin{proof}
If $\Pi_{I} = 0$, then $G(\Pi_{I}) = O(n)$. So, assume $\Pi_{I} \neq 0$. By the same argument as in Theorem \ref{t:3-structure}, Case 3, intersecting $G(\Pi_I)$ with any of the subgroups $(3.a)$, $(3.b)$ and $(3.c)$ implies that, eventually, the resulting symmetry group must be of one the forms $(3.a)$, $(3.b)$ or $(3.c)$. Therefore, we may only look into the intersection of $G(\Pi_I)$ with subgroups of the form $(3.d)$. Assume that the latter are expressed with respect to an orthonormal basis $r_1, r_2, r_3$.

Since $\Pi_{I} \in so(3)$, then there exists $S_{I}:=(p_1, p_2, p_3) \in O(3)$ such that
$\Pi_{I} = S_{I} \textnormal{diag}(L_{s},0)S^{*}_{I}$, where $L_{s} \neq 0$ has the form (\ref{e:L_s}).
Therefore, by Theorem \ref{t:comm_general}, we have
\begin{equation}\label{e:G(Pi_I)_types-3}
G(\Pi_I) = S_{I} \textnormal{diag}(SO(2),\pm 1)S^*_{I}.
\end{equation}
Note that, for a matrix of the form $O \textnormal{diag}(SO(2),\pm
1)O^* = OU\textnormal{diag}(e^{i \theta},e^{-i \theta},\pm 1)U^*O^*$,
$\theta \in (0,2 \pi) \backslash \{\pi\}$, where $U = \textnormal{diag}(U_2,1)$ and $U_2$ is as in (\ref{e:U2}), only the eigenvalue 1 (or $-1$) is associated
with a purely real eigenvector. Thus, we can break up the rest of the proof into two cases.
Assume that $\textnormal{span}_{\bbR}\{p_3\}=\textnormal{span}_{\bbR}\{r_3\}$. Then we can write $G(\Pi_I) = (r_1,r_2,r_3)\textnormal{diag}(SO(2),\pm 1)(r_1,r_2,r_3)^*$.
Thus, $G(\Pi_I) \cap (r_1,r_2,r_3)\textnormal{diag}(O(2),\pm 1)(r_1,r_2,r_3)^* = G(\Pi_I)$, which gives $(3.f)$. Alternatively, assume $\textnormal{span}_{\bbR}\{p_3\} \neq \textnormal{span}_{\bbR}\{r_3\}$. In this case, one can argue exactly as in the proof of Theorem \ref{t:3-structure}, Case 3, $(i)$-$(iii)$, to obtain the same sets of candidates ${\mathcal C}$. Thus, the claim follows. $\Box$\\
\end{proof}

\begin{example}\label{ex:3-3.f} (Type $(3.f)$)
Analogously to Example \ref{ex:2-rotational}, consider an OFBM with parameters
$$
D = d I, \quad \Re(AA^*) = I, \quad\Im(AA^*) = \textnormal{diag}(L,0), \quad L \in so(2) \backslash\{0\},
$$
where $d$ is real. Then, $G(\Pi_x) = O(3)$ and $G(\Pi_{I})$ is as in (\ref{e:G(Pi_I)_types-3}).
\end{example}

Theorems \ref{t:3-structure} and \ref{t:3-structure_nontimerevers} stand in contrast with Theorem \ref{t:2-structure} in that they show the much greater wealth of possible symmetry groups in dimension 3 as compared to dimension 2. In a certain sense, this enhances the claim of Theorem \ref{t:minimal_is_topol_general} in that, notwithstanding the increasing complexity of the possible symmetry structures as dimension increases, minimal type symmetry groups remain the topologically general case for \textit{any} dimension.

We now provide the tangent spaces and exponent sets for each
symmetry group with non-trivial tangent space. The proof is along
the lines of that for Proposition \ref{p:2-exponents}.
\begin{proposition}
Under the assumptions of Theorem \ref{t:3-structure_nontimerevers}, for the
symmetry groups associated with non-trivial tangent spaces, the
tangent spaces, commuting exponents $H_0$ and sets of exponents
have the form:
\begin{itemize}
\item [(3.d)] for some orthonormal $p_1$, $p_2$, $p_3$, and the
associated matrix $S :=(p_1, \hspace{1mm} p_2,\hspace{1mm} p_3)$,
$$
T(G_H) = WS \textnormal{diag}(so(2),0)S^*W^{-1}, \quad H_0 = WSU\textnormal{diag}(h_1,\overline{h}_1,h_2)U^*S^*W^{-1},
$$
$$
{\mathcal E}(B_H) =WSU(\textnormal{diag}(h_1,\overline{h}_1,h_2)+
\textnormal{diag}(so(2),0))U^* S^*W^{-1},
$$
where $U = \textnormal{diag}(U_2,1)$ and $U_2$ is as in
(\ref{e:U2}), and $h_1 \in \bbC$, $h_2 \in \bbR$;

\item [(3.e)] $T(G_H) = T(SO(3)) = Wso(3)W^{-1}$, $H_0 = h_0 I$, ${\mathcal E}(B_H) = h_0 I + W so(3) W^{-1}$;

\item [(3.f)] the same as for (3.d).
\end{itemize}
\end{proposition}
\begin{proof}
For type $(3.d)$, just note that $T(G_H) =
T(\textnormal{Rot}_{\theta}(p_3))=WS\textnormal{diag}(so(2),0)S^*W^{-1}$, from
which $H_0$ and ${\mathcal E}(B_{H})$ promptly follow. The same
argument holds for type $(3.f)$.

The case of type $(3.e)$ is straightforward, since $T(G_H) =
T(SO(3))$. $\Box$\\
\end{proof}

\begin{remark}
In general dimension $n$, there are no additional difficulties in
describing the structure of groups $G(\Pi)$ for a \textit{fixed}
symmetric matrix $\Pi$. Equivalently, one can generalize Proposition \ref{p:3-G(Pi)}
to the context of dimension $n$ without much effort. Nevertheless,
it is cumbersome to describe the structure
of intersections $G(\Pi_1)\cap G(\Pi_2)$, which is needed for the full
characterization of symmetry
groups $G_H$ as in (\ref{e:symmetry-group-OFBM-main}) and
(\ref{e:symmetry-group-OFBM-main-2}). At this point, a full
description of symmetry groups in general dimension $n$ is an open
question.
\end{remark}

\begin{remark}\label{r:SO(n)_not_a_symmgroup}
The classification given in Theorem \ref{t:2-structure} stands in contrast
with the fact that $SO(n)$ cannot be a symmetry group for
$\bbR^n$-valued random vectors (Billingsley \cite{billingsley:1966}). In particular, $SO(2)$ is not a maximal element of its equivalence class of subgroups in the sense of Meerschaert and Veeh \cite{meerschaert:veeh:1995}, p.\ 2 (not to be confused with the symmetry group of maximal type in Theorem \ref{t:2-structure}). However, it turns out that Billingsley's result is actually \textit{almost}
true for OFBMs, and more generally, proper zero mean Gaussian processes. In other words, for the latter class of processes, $SO(n)$ can only be a symmetry group when $n=2$ (cf.\ Theorem \ref{t:3-structure_nontimerevers}). Indeed, without loss of generality, assume $W=I$. Then, it suffices to show that $SO(n) \subseteq G(X)$ implies that $O(n) = G(X)$ when $n \geq 3$. However, the latter equivalence is a consequence of Proposition \ref{p:Centr=O(n)=>lambdaI} in the appendix.
\end{remark}

\section{On integral representations of OFBMs with multiple
exponents} \label{s:integral-repres}

In this section, we show that when an OFBM has multiple exponents,
the matrix $A$ in (\ref{e:spectral-repres-OFBM}) can be chosen the same, no matter
what matrix exponent is used in the parametrization. We also show that, by contrast, such invariance of
the parametrization does \textit{not} hold for the so-called time
domain representation of OFBM.

We first consider the latter point. Under (\ref{e:eigen-assumption}) and $\Re(h) \neq 1/2$
for any eigenvalue $h$ of $H$, the OFBM $\{B_{H}(t)\}_{t \in \bbR}$
also admits an integral representation in the time domain, i.e.,
\begin{equation}\label{e:integ_repres_timedom}
\{B_{H}(t)\}_{t \in \bbR} \stackrel{{\mathcal L}}= \Big\{
\int_{\bbR} (((t-u)^{H - \frac{1}{2}I}_{+} - (-u)^{H -
\frac{1}{2}I}_{+})M_{+} + ((t-u)^{H - \frac{1}{2}I}_{-} - (-u)^{H
- \frac{1}{2}I}_{-})M_{-})B(du)\Big\}_{t \in \bbR},
\end{equation}
where $M_{+}$, $M_{-} \in M(n,\bbR)$, and $\{B(u)\}_{u \in \bbR}$ is a vector-valued
process consisting of independent Brownian motions and such that
$EB(du)B(du)^* = du$ (Didier and Pipiras \cite{didier:pipiras:2009}). The following example
shows that, in general, the matrix parameters $M_{+}$, $M_{-}$ cannot be chosen
independently of the exponent.

\begin{example}\label{ex:timedom_params_affected_symm}
Consider a bivariate OFBM $B_{H}$ with the time domain representation
 (\ref{e:integ_repres_timedom}), where $D = dI$, $d \in (-1/2,1/2) \backslash \{0\}$ (or $H=hI$, $h \in (0,1)\backslash\{1/2\}$),
 $M_{+} = O \in SO(2)$ and $M_{-} = I$. Since rotation matrices
 commute, it follows directly from (\ref{e:integ_repres_timedom}) that
 \begin{equation}\label{e:SO(2)_subseteq_GH}
 SO(2) \subseteq G_{H}.
 \end{equation}
The relation (\ref{e:SO(2)_subseteq_GH}) implies that $T(SO(2)) = so(2) \subseteq T(G_H)$.
Hence, in view of (\ref{e:exponents}),
\begin{equation}\label{e:H+Lc}
H + L_c, \quad c \in \bbR,
\end{equation}
are the exponents of the OFBM $B_H$, where $L_c \in so(2)$ is given in (\ref{e:L_s}). Thus,
the OFBM $B_H$ has the time domain representation
\begin{equation}\label{e:integ_repres_timedom_Lc}
\{B_{H}(t)\}_{t \in \bbR} \stackrel{{\mathcal L}}= \Big\{ \int_{\bbR} (((t-u)^{D+L_c}_{+} - (-u)^{D+L_c}_{+})M_{+} +
((t-u)^{D+L_c}_{-} - (-u)^{D+L_c}_{-})M_{-})B(du) \Big\}_{t \in \bbR},
\end{equation}
where $M_{+}=M_{+}(c)$, $M_{-}=M_{-}(c)$. We want to show that one cannot generally take
the original parameters $M_{+}=O$, $M_{-}=I$ in the representation (\ref{e:integ_repres_timedom_Lc}).

Arguing by contradiction, suppose that $M_{+}=O$, $M_{-}=I$ in
(\ref{e:integ_repres_timedom_Lc}) lead to the same OFBM for any $c
\in \bbR$. In the spectral domain, these processes have the
representation
\begin{equation}\label{e:integ_repres_specdom_Lc}
\int_{\bbR}\frac{e^{itx}-1}{ix}(x^{-D - L_c}_{+}A  + x^{-D - L_c}_{-}\overline{A})\widetilde{B}(dx),
\end{equation}
where $\widetilde{B}(dx)$ is as in (\ref{e:spectral-repres-OFBM}), and
$$
A = \frac{1}{\sqrt{2\pi}}\Gamma(D+L_c+I)(e^{-i \pi (D+L_c)/2}O +
e^{i \pi (D+L_c)/2}I)
$$
(see Theorem 3.2 and its proof in Didier and Pipiras \cite{didier:pipiras:2009}). Note that
$x^{-L_c}$ commutes with $A$, $\overline{A}$ and hence $(-L_c)$ can be removed from the exponents
of $x_+$, $x_-$ in (\ref{e:integ_repres_specdom_Lc}). Then, if (\ref{e:integ_repres_specdom_Lc})
is the same process for all $c \in \bbR$, the matrix
\begin{eqnarray}\label{e:AA*_when_timedom_notunique}
(2 \pi)AA^* & =&  \Gamma(dI + L_{c}+I) (e^{-i \pi (dI + L_c)/2}O + e^{i \pi (dI + L_c)/2}I)\cdot \nonumber \\
& & \hspace{1mm}\cdot(O^* e^{i \pi (dI + L^*_c)/2} + I e^{-i \pi (dI + L^*_c)/2})  \Gamma(dI + L_{c}+I)^{*} \nonumber  \\
& = & \Gamma(dI + L_{c}+I)\Gamma(dI + L_{c}+I)^{*}(e^{- i \pi dI}O + e^{i \pi dI}O^* + e^{-i \pi L_c} + e^{i \pi L_{c}} )
\end{eqnarray}
does not depend on $c$. Note that
$$
\Gamma(dI + L_c+I)\Gamma(dI + L_c+I)^* = U_2
\textnormal{diag}(|\Gamma(d+ic +1)|^2,|\Gamma(d+ic +1)|^2)U^*_{2} =
|\Gamma(d+ic +1)|^2I,
$$
where $U_2$ is as in (\ref{e:U2}) and $\Gamma(d+ic+1)$ is the univariate Gamma function evaluated at
$d+ic +1 \in \bbC$. Writing $O = U_2 \textnormal{diag}(e^{i \beta},e^{-i \beta})U^*_2$, for some
$\beta \in (0, 2 \pi)\backslash\{\pi\}$, the matrix (\ref{e:AA*_when_timedom_notunique}) becomes
$U_2 \textnormal{diag}(f(d,c,\beta),\overline{f(d,c,\beta)})U^*_2$, where
$$
f(d,c,\beta):=|\Gamma(d+ic +1)|^2\Big(e^{-i
\pi d}e^{i\beta}+e^{i \pi d}e^{-i\beta}+ e^{\pi c}+e^{-\pi
c}\Big).
$$
However, the function $f(d,c,\beta)$ does depend on $c$, as can be easily verified (contradiction).
\end{example}

The following result shows that, for a given OFBM, one can take the same parameter $A$
in the spectral representation (\ref{e:spectral-repres-OFBM}) for all exponents $H$ of the OFBM
in question.

\begin{theorem}\label{t:AA*_is_unique}
Let $B_H$ be an OFBM having the spectral representation
(\ref{e:spectral-repres-OFBM}). If $H_{\lambda}, H_{\eta} \in
{\mathcal E}(B_H)$ and $A_{\lambda}$, $A_{\eta}$ are the two
matrix parameters in (\ref{e:spectral-repres-OFBM}) associated
with $H_{\lambda}$, $H_{\eta}$, respectively, then
\begin{equation}\label{e:AA*_irrespective_of_H}
A_{\lambda}A^*_{\lambda} = A_{\eta}A^*_{\eta}.
\end{equation}
In particular, one may choose the same matrix parameter $A$ in (\ref{e:spectral-repres-OFBM}) for
every choice of $H \in {\mathcal E}(B_H)$.
\end{theorem}
\begin{proof}
It is enough to show (\ref{e:AA*_irrespective_of_H}) with a commuting exponent $H_{\eta}:=H_0$ (see (\ref{e:H0_def})) and the associated matrix $A_{\eta}:=A_0$. For simplicity,
let $H=H_{\lambda}$, $A=A_{\lambda}$. We know that
$$
H-H_0 = D - D_0 =: \Delta \in W {\mathcal L}_{0}W^{-1} = T(G_H),
$$
where ${\mathcal L}_{0} \subseteq so(n)$. We can thus write $\Delta = WLW^{-1}$ with
$L \in so(n)$. The uniqueness of the spectral density of OFBM implies that, for $x > 0$,
$$
x^{-D}AA^*x^{-D^*}=x^{-D_0}A_0 A^*_0 x^{-D^*_0},
$$
or
$$
x^{-(D_0+\Delta)}AA^*x^{-(D_0+\Delta)^*}=x^{-D_0}A_0 A^*_0 x^{-D^*_0}.
$$
Since $D_0$ is a commuting exponent and $\Delta \in T(G_H)$, then $D_0$ and $\Delta$ commute.
Hence,
$$
x^{-D_0}x^{-\Delta}AA^* x^{-\Delta^*}x^{-D^*_{0}} = x^{-D_0}A_0 A^*_0 x^{-D^*_{0}}
$$
and $x^{-\Delta}AA^*x^{-\Delta^*} = A_0 A^{*}_{0}$, i.e., $x^{-L}W^{-1}AA^*W^{-1}x^{L}= W^{-1}A_0 A^*_0 W^{-1}$. By differentiating with respect to $x$, we 
further obtain $L(W^{-1}AA^*W^{-1}) = (W^{-1}AA^*W^{-1})L$, that is, $L$ and $W^{-1}AA^*W^{-1}$ commute. Then
$W^{-1}AA^* W^{-1} = W^{-1}A_0 A^*_{0}W^{-1}$ or $AA^* = A_0A^*_0$. The last statement of the theorem follows from (\ref{e:AA*_irrespective_of_H}). $\Box$\\
\end{proof}

\appendix

\section{Auxiliary results on matrix commutativity}
\label{s:auxiliary}
%\section{On the commutativity of matrices} \label{a:commutativity}
We begin by proving Corollary \ref{c:commute=>invariant} and Lemma \ref{l:irreducible}. The argument draws upon Theorem \ref{t:comm_general}.

\bigskip {\sc Proof of Corollary \ref{c:commute=>invariant}}:
Without loss of generality, we can assume that $j_1 = 1, \hdots, j_k = k$. Since $M$ commutes with $A$, by Theorem \ref{t:comm_general} we can write $M = O K O^*$, where $K = \textnormal{diag}(K_{11},K_{22})$ and $K_{11} \in M(k,\bbR), K_{22} \in M(n-k,\bbR)$. The rest of the claim immediately follows. $\Box$\\

\bigskip {\sc Proof of Lemma \ref{l:irreducible}}:
Consider the spectral decomposition $A = O \Upsilon O^*$, $O = (o_1,\hdots,o_{n}) \in O(n)$ and $\Upsilon $ is diagonal. Then, since $M$ and $A$ commute, by Corollary \ref{c:commute=>invariant} and the assumption that $A \in {\mathcal S}_{\neq}$ we have $M = O\textnormal{diag}(\lambda_1,\hdots,\lambda_n)O^*$ for some $\lambda_1,\hdots,\lambda_n \in \bbR$. In particular, $M \in {\mathcal S}(n,\bbR)$. So, assume by contradiction that $\lambda_{j_1}=\hdots = \lambda_{j_k}$ (possibly $k = 1$), where $\lambda_{j_1} \neq \lambda_i$ for any other eigenvalue $\lambda_i$ of $M$. Without loss of generality, we can assume that $j_1 = 1, \hdots, j_k = k$. But since $B$ commutes with $M$, again by Corollary \ref{c:commute=>invariant} $\textnormal{span}_{\bbR}\{o_{1},\hdots, o_{k}\}$ is a $k$-dimensional invariant subspace of $B$, which contradicts the assumption. Therefore, $\lambda_1= \hdots =\lambda_n$, as claimed. $\Box$\\

The next proposition is used in Theorem \ref{t:maximal} and Remark \ref{r:SO(n)_not_a_symmgroup}. It shows that, for $n \geq 3$, the group $SO(n)$ is so rich
that only a matrix which is a multiple of the identity can contain it in its centralizer. For $n=2$, one needs
to consider instead the entire orthogonal group.

\begin{proposition}\label{p:Centr=O(n)=>lambdaI}
Let $\Gamma \in M(n,\bbR)$. Then $\Gamma = \lambda \hspace{0.5mm}I$, $\lambda \in \bbR$,
if one of the following assumptions holds:
\begin{enumerate}
\item [(i)] for $n=2$, if ${\mathcal C}(\Gamma) \supseteq O(n)$;
\item [(ii)] for $n \geq 3$, if ${\mathcal C}(\Gamma) \supseteq
SO(n)$.
\end{enumerate}
\end{proposition}
\begin{proof}
We only prove $(ii)$.

%It is a consequence of slight modifications of Lemma \ref{l:irreducible} and Corollaries \ref{c:invar=>zeroes_in_S_and_L} and \ref{c:invar_subs=>eigenvec}. We restate them appropriately, but omit the proofs since they can be worked out with little effort. So, let ${\mathcal N}(n,\bbF)$ be the set of normal matrices, $\bbF = \bbR$ or $\bbC$, and let ${\mathcal N}_{\neq}(n,\bbF)$ be the subset of normal matrices with pairwise distinct eigenvalues. Note that the next statements involve complex invariant subspaces.
%
%\begin{itemize}
%\item Let $A,B \in {\mathcal N}_{\neq}(n,\bbR)$. Assume $A$, $B$ have no $k$-dimensional invariant complex subspaces in common, $k = 1,\hdots,n-1$. If $M \in M(n,\bbR)$ commutes with both $A$ and $B$, then $M$ is a scalar matrix.
%\item Let $M \in {\mathcal N}(n,\bbR)$. If the orthonormal vectors $u_1,\hdots,u_n$ generate an invariant subspace of $M$, then
%$$
%M = U
%\left(\begin{array}{cc}
%M_{11} & 0\\
%0 & M_{22}
%\end{array}\right)
%U^*,
%$$
%where $U \in U(n)$, $M_{11} \in {\mathcal N}(k,\bbC)$, $M_{22} \in {\mathcal N}(n-k,\bbC)$.
%\item Let $M \in {\mathcal N}(n,\bbR)$. Any $k$-dimensional complex invariant subspace of $M$ is generated by $k$ eigenvectors of $M$. Moreover, its perp is also a $(n-k)$-dimensional invariant complex subspace.
%\end{itemize}

First, let $n$ be odd. We have that $\Gamma$ commutes with all rotation matrices of the form
$$
\textnormal{diag}(-1,\hdots,-1,1,-1,\hdots,-1),
$$
(i.e., the eigenvalue 1 is taken in every ``position", all the remaining ones being filled in by $-1$; there is an even number of entries with $-1$, so the determinant is 1). By Corollary \ref{c:commute=>invariant}, this implies that the Euclidean vectors $e_1, \hdots, e_n$  are all eigenvectors of $\Gamma$. Consequently, $\Gamma$ is diagonalizable and we can write $\Gamma = \textnormal{diag}(\lambda_1,\hdots,\lambda_n)$. Since $\Gamma \in M(n,\bbR)$, then $\lambda_1,\hdots, \lambda_n \in \bbR$.

Furthermore, more generally we have that $\Gamma$ commutes with all orthogonal matrices of the form
$$
P \textnormal{diag}(-1,\hdots,-1,1,-1,\hdots,-1) P^*, \quad P \in O(n).
$$
By the same reasoning as above, this implies that all $p \in S^{n-1}$ are eigenvectors of $\Gamma$. Consequently, every vector $v \in \bbR^n$ is an eigenvector of $\Gamma$.

We claim that $\lambda_1 = \hdots = \lambda_n$. Assume by contradiction that this does not hold. Without loss of generality, we can take $\lambda_1 \neq \lambda_2$. Then
$\Gamma(e_1 + e_2) = \lambda_1 e_1 + \lambda_2 e_2$, but also $\Gamma(e_1 + e_2) = \gamma (e_1 + e_2)$ for some $\gamma \in \bbR$, since $e_1 + e_2$ is an eigenvector of $\Gamma$ and $\Gamma \in M(n,\bbR)$. Thus,
$(\gamma - \lambda_1)e_1 = (\lambda_2 - \gamma)e_2$. Since $e_1$ and $e_2$ are linearly independent, $\gamma = \lambda_1$ and $\lambda_2 = \gamma$, which contradicts the assumption that $\lambda_1 \neq \lambda_2$.

Now, let $n$ be even. For notational simplicity, take $n = 4$. $\Gamma$ commutes with $\textnormal{diag}(1,1,-1,-1) \in SO(4)$. Thus, by Theorem 2.1, $\Gamma$ has the form
$$
\Gamma = \left(\begin{array}{cccc}
a & b & & \\
c & d & & \\
  &   & e & f \\
  &   & g & h \\
\end{array}\right),
$$
all entries being in $\bbR$. Likewise, $\Gamma$ commutes with $\textnormal{diag}(1,-1,-1,1) \in SO(4)$. Thus, by Theorem 2.1, $\Gamma$ has the form
$$
\Gamma = \left(\begin{array}{cccc}
i &  & & j \\
  & m  &  n & \\
  & o  &  p  &  \\
k   &   &  & l \\
\end{array}\right),
$$
all entries being in $\bbR$ again. Consequently, $\Gamma$ must have the form $\Gamma = \textnormal{diag}(a,d,e,h)$. In other words, $e_1, \hdots, e_4$ 
are eigenvectors of $\Gamma$ (and $\Gamma$ is diagonalizable). Now, note that $\Gamma$ commutes with $P \textnormal{diag}(1,1,-1,-1) P^*, P \textnormal{diag}(1,-1,-1,1) P^* \in SO(4)$ for all $P \in O(4)$. Thus, we conclude that every $p \in S^{3}$ is an eigenvector of $\Gamma$, and thus that every $v \in \bbR^4$ is an eigenvector of $\Gamma$. Consequently, by the same reasoning as for the last argument for the case of odd $n$, $\Gamma = \lambda I$ for some $\lambda \in \bbR$.

This argument can be extended to the case of general even $n$ by noting that $\Gamma$ commutes with all rotation matrices of the form
$$
\textnormal{diag}(-1,\hdots,-1,1,1,-1,\hdots,-1),
$$
where the ``window" associated with the pair (1,1) ``slides" all the way down the main diagonal. Thus, by again considering all conjugacies $P \in O(n)$, we can see that the same argument applies.

$\Box$\\
\end{proof}

\begin{remark}
The claim can be extended to the case of $\Gamma \in M(n,\bbC)$ (which is used in Didier and Pipiras \cite{didier:pipiras:2009}, Proposition 5.3). When $n$ is odd, we obtain instead that $\Gamma = \textnormal{diag}(\lambda_1,\hdots,\lambda_n)$, where $\lambda_1,\hdots,\lambda_n \in \bbC$. Moreover, we conclude that every $v \in \bbR^n$ (not $\bbC^n$) is an eigenvector. Thus, if $\lambda_1 \neq \lambda_2$, we have that
$$
\gamma (e_1 + e_2 ) = \Gamma(e_1 + e_2) = \lambda_1 e_1 + \lambda_2 e_2
$$
for some $\gamma \in \bbC $, which is a contradiction. The proof for the case of even $n$ can be obtained by similar adaptation of the argument for the real case.
\end{remark}

As a consequence of Proposition \ref{p:Centr=O(n)=>lambdaI}, we obtain the next lemma, which is used in Theorem \ref{t:maximal}.

\begin{lemma}\label{l:WO(n)W^(-1)=W2OW^(-1)2}
If
$$
W_1 O(n) W^{-1}_1 = W_2 {\mathcal O} W^{-1}_2,
$$
where $W_1$ and $W_2$ are positive definite matrices, and ${\mathcal O}$ is a subset (subgroup) of $O(n)$, then ${\mathcal O} =  O(n)$.
\end{lemma}
\begin{proof}
Write $Q O(n) Q^{-1} =  {\mathcal O}$, where $Q := W^{-1}_2 W_1$. Thus, for any $O_1\in O(n)$, $O_2 := Q O_1 Q^{-1} \in O(n)$, i.e.,
$(Q O_1 Q^{-1})(Q O_1 Q^{-1})^{*} = I$. Thus $Q O_1 Q^{-1} = (Q^*)^{-1} O_1 Q^*$, or, equivalently,
$(Q^* Q) O_1 = O_1 (Q^* Q)$. Thus, by Proposition \ref{p:Centr=O(n)=>lambdaI}, there exists $\lambda > 0$ such that $Q^* Q = \lambda I$. Therefore,
$W_2 = W_1 \lambda^{-1/2}$, from which the claim follows. $\Box$\\
\end{proof}

The next lemma is used in Theorem \ref{t:3-structure}.

\begin{lemma}\label{l:third_eigenvec_of_O}
Let $O \in O(3)$. Assume $O$ has two (real) linearly independent eigenvectors $p_1$, $p_2$, both of which are associated with real eigenvalues. Then $O$ has a third (real) eigenvector $p_3$ such that $p_3 \perp \textnormal{span}_{\bbR}\{p_1,p_2\}$.
\end{lemma}
\begin{proof}
It is clear that, under the assumptions, $O$ has only real eigenvalues. We break up the argument into subcases.

\begin{itemize}
\item [(i)] $p_1$, $p_2$ are eigenvectors associated with the same real eigenvalue (without loss of generality, assume the latter is 1): if the third eigenvalue of $O$ is 1, then the proof ends. Otherwise, since $O \in O(3)$, eigenvectors associated with different eigenvalues are orthogonal. Thus, there exists a vector $p_3$ (associated with $-1$) which is orthogonal to $\textnormal{span}_{\bbR}\{p_1,p_2\}$. Moreover, since the eigenvalue in question is real and $O \in O(3)$, one can assume that $p_3$ is real.
\item [(ii)] $p_1$, $p_2$ are eigenvectors associated with different real eigenvalues (without loss of generality, assume such eigenvalues are 1 and $-1$, respectively): since $O \in O(3)$, $p_1 \perp p_2$. Moreover, since $O \in O(3)$, then the eigenspace associated with one of the eigenvalues 1 or $-1$ has dimension 2. Assume without loss of generality that such eigenvalue is 1. Thus, there must exist another eigenvector $p_3$ associated with 1. Moreover, since $O \in O(3)$, we can assume that $p_3$ is real. We obtain that $p_3 \perp p_2$. Also, it is clear that one can choose $p_3$ so that $p_3 \perp p_1$, and thus $p_3 \perp \textnormal{span}_{\bbR}\{p_1,p_2\}$. 
\end{itemize}
$\Box$
\end{proof}

\section{On the convergence of eigenvalues and eigenvectors}

The next two lemmas are used in Lemmas \ref{e:open-dense_separate} and \ref{l:graph(L)_is_open}.

The first one, whose proof is straightforward, establishes the sense in which the convergence of matrices implies the convergence of
the eigenvalues.

The second one shows that, under more stringent assumptions, the convergence of matrices also implies the convergence of the eigenvectors (in the
specific sense described in the claim).

\begin{lemma}\label{l:conv_of_eigenvalues}
Let $\{A_k\}_{k \in \bbN}, A_0 \in M(n,\bbC)$, and assume that $A_k \rightarrow A_0$. Then, the eigenvalues of $A_k$ converge to those of $A_0$, i.e.,
one can form a sequence $\{(\lambda^{1}_{k},\hdots,\lambda^{n}_{k})\}_{k \in \bbN} \subseteq \bbC^n$ of eigenvalues of $A_k$, $k \in \bbN$, such that
\begin{equation}\label{e:conv_of_eigenvalues}
(\lambda^{1}_{k},\hdots,\lambda^{n}_{k}) \rightarrow (\lambda^{1}_{0},\hdots,\lambda^{n}_{0}),
\end{equation}
where the vector on the right-hand side of (\ref{e:conv_of_eigenvalues}) consists of eigenvalues of $A_0$.
\end{lemma}

\begin{lemma}\label{l:Pk_conv_to_P}
Let $\{A_k\}_{k \in \bbN}, A_0 \in M(n,\bbC)$ be such that $A_0$ has pairwise distinct eigenvalues. Assume that $A_k \rightarrow A_0$. Then, for the sequence of eigenvalues in (\ref{e:conv_of_eigenvalues}), there exists a sequence of conjugacies $\{P_k\}_{k \in \bbN} \subseteq GL(n,\bbC)$ such that, for large $k$,
$$
A_{k} = P_k \textnormal{diag}(\lambda^1_{k},\hdots,\lambda^n_{k}) P^{-1}_{k}
$$
and $P_k \rightarrow P$ for some $P \in GL(n,\bbC)$. Moreover, the columns of the limiting matrix $P$ are, in fact, eigenvectors of $A_0$.
\end{lemma}
\begin{proof}
By Lemma \ref{l:conv_of_eigenvalues}, the eigenvalues of $A_k$ converge to those of $A_0$. Consequently, for large enough $k$, the eigenvalues of $A_k$, which we can write as $\lambda^{j}_{k}, j=1,\hdots, n$, are pairwise distinct. For all $j=1,\hdots,n$, choose a basis $p^{1}_0,\hdots,p^{n}_0 \in S^{2n-1}$ of eigenvectors of $A_0$.

Now for large enough $k$ and for each $j$, let $V_{j,k} := \{p \in \bbC^n: A_k p = \lambda^{j}_{k} p  \}$ be the one-dimensional eigenspace associated with the eigenvalue $\lambda^{j}_{k}$. $V_{j,k}$ is well-defined for large $k$ because $\lambda^{1}_{k},\hdots, \lambda^{n}_{k}$ are pairwise distinct. Now let $p^{j}_k$ be one eigenvector in $S^{2n-1}$ that is closest to $p^{j}_0$, i.e.,
\begin{equation}\label{e:contra_pjk_is_eigenvec_argmin}
p^{j}_k \in \textnormal{argmin}_{p \in V_{j,k} \cap S^{2n-1}} \norm{p - p^{j}_0 }.
\end{equation}
The compactness of $V_{j,k} \cap S^{2n-1}$ implies that a solution must exist.

Fix any $j = 1,\hdots,n$. Assume that $\{p^{j}_k\}$ does not converge to $p^j_0$, i.e., there exists $\varepsilon_0 > 0$ and a subsequence $\{p^{j}_{k'}\}$ such that
\begin{equation}\label{e:contra_argmin_eigenvec_not_conv}
\norm{p^{j}_{k'} - p^j_0} \geq \varepsilon_0.
\end{equation}
Since $\{p^{j}_{k'}\} \subseteq S^{2n-1}$, one can extract a further subsequence $\{p^{j}_{k''}\}$ which is convergent, i.e.,
\begin{equation}\label{e:contra_subseq_conv_p''}
p^{j}_{k''} \rightarrow p'' \in S^{2n-1}.
\end{equation}
Therefore, $A_{k''}p^{j}_{k''} = \lambda^{j}_{k''} p^{j}_{k''}$, where $A_{k''}p^{j}_{k''} \rightarrow A_0 p'', \quad \lambda^{j}_{k''} p^{j}_{k''}\rightarrow \lambda^{j}_0 p''$.
Thus, $p'' \in S^{2n-1}$ is an eigenvector of $A_0$ associated with the eigenvalue $\lambda^{j}_0$. Then, since the eigenvalues of $A_0$ are pairwise distinct, we can write $p'' = e^{i \theta} p^j_0$ for some $\theta \in [-\pi,\pi)$. Note that $e^{-i\theta}p^{j}_{k''}$ is an eigenvector of $A_{k}$ associated with the eigenvalue $\lambda^{j}_k$. Thus, from (\ref{e:contra_pjk_is_eigenvec_argmin}), (\ref{e:contra_argmin_eigenvec_not_conv}) and (\ref{e:contra_subseq_conv_p''}), we obtain
$$
0 < \varepsilon_0 \leq \norm{p^{j}_{k''} - p^{j}_{0}} \leq \norm{e^{-i \theta} p^{j}_{k''} - p^{j}_{0}} = |e^{-i \theta}|\norm{ p^{j}_{k''} - e^{i\theta} p^{j}_{0}} = \norm{ p^{j}_{k''} - p''} \rightarrow 0, \quad k'' \rightarrow \infty
$$
(contradiction). Therefore, $p^{j}_{k} \rightarrow p^{j}_{0}$, $j = 1,\hdots,n$, and we can define $P_{k} = (p^{1}_{k},\hdots,p^{n}_k)$. Since such eigenvectors must be linearly independent, we can write $A_k = P_k \textnormal{diag}(\lambda^{1}_k,\hdots,\lambda^{n}_{k})P^{-1}_k$, and $P_k \rightarrow P$, where $P = (p^{1}_0, \hdots, p^{n}_0)$. $\Box$\\
\end{proof}

\section{Additional results on matrix representations}

The results in this section are used in Section \ref{s:min-max}. The proof of Lemma \ref{l:invar=>zeroes} is omitted because the argument is standard.

\begin{lemma}\label{l:invar=>zeroes}
Let $M$ be a matrix in $M(n,\bbR)$ and let $p_1,\hdots,p_k$ be linearly independent vectors in $\bbR^n$. Assume $\textnormal{span}_{\bbR}\{p_1,\hdots,p_k\}$ is an invariant subspace of $M$. For $p_{k+1},\hdots,p_n$ such that $p_1,\hdots,p_k, p_{k+1},\hdots,p_n $ is a basis of $\bbR^n$, define the conjugacy matrix $P = (p_1,\hdots,p_k,p_{k+1},\hdots,p_n) \in GL(n,\bbR)$. Then we obtain the representation
$$
M = P \left(\begin{array}{cc}
M_{11} & M_{12} \\
{\mathbf 0}  & M_{22}
\end{array}\right)P^{-1},
$$
where $M_{11} \in M(k,\bbR), M_{22} \in M(n-k,\bbR), M_{12} \in M(k,n-k,\bbR)$.
\end{lemma}

\begin{corollary}\label{c:invar=>zeroes_in_S_and_L}
Let $S \in {\mathcal S}(n,\bbR)$. If the real orthonormal vectors $o_1,\hdots,o_k$ generate an invariant subspace of $S$, then
\begin{equation}\label{e:S_repres_invar_subspace}
S = O \left(\begin{array}{cc}
S_{11} & {\mathbf 0} \\
{\mathbf 0}  & S_{22}
\end{array}\right)O^{*},
\end{equation}
where $O := (o_1,\hdots,o_k,o_{k+1},\hdots,o_n) \in O(n)$, $S_{11} \in {\mathcal S}(k,\bbR), S_{22} \in {\mathcal S}(n-k,\bbR)$. An analogous conclusion holds if we replace $S \in {\mathcal S}(n,\bbR)$ with $L \in so(n)$.
\end{corollary}

\begin{corollary}\label{c:invar_subs=>eigenvec}
Let $S \in {\mathcal S}(n,\bbR)$. Each $k$-dimensional real invariant subspace of $S$ is generated by a set of $k$ real eigenvectors of $S$. Moreover, the orthogonal subspace is also a $(n-k)$-dimensional invariant subspace.
\end{corollary}
\begin{proof}
Assume that $\textnormal{span}_{\bbR}\{o_1,\hdots,o_{k}\}$ is a $k$-dimensional real invariant subspace of $S$, where, by the Gram-Schmidt algorithm, we can assume that the real vectors $o_1,\hdots,o_k$ are orthonormal. Then the representation (\ref{e:S_repres_invar_subspace}) holds for $S$, where the first $k$ vectors of $O$ are $o_1,\hdots, o_k$. Consider the spectral decompositions $S_{11}= P_{11} \textnormal{diag}(\lambda_1,\hdots,\lambda_k) P^*_{11}$ and $S_{22}=P_{22} \textnormal{diag}(\lambda_{k+1},\hdots,\lambda_{n}) P^*_{22}$. We can write
$$
S = O \textnormal{diag}(P_{11},P_{22}) \textnormal{diag}(\lambda_1,\hdots,\lambda_k,\lambda_{k+1},\hdots,\lambda_{n})  \textnormal{diag}(P^{*}_{11},P^*_{22}) O^*.
$$
Thus, the column vectors in the matrix $(o_1,\hdots,o_k)P_{11}$ are (linearly independent) eigenvectors of $S$. Moreover, since they are all linear combinations of the vectors $o_1, \hdots, o_k$, 
\begin{equation}\label{e:span_orthogonal_basis}
\textnormal{span}_{\bbR}\{(o_1,\hdots,o_k)P_{11}\} \subseteq \textnormal{span}_{\bbR}\{o_1,\hdots,o_k\} .
\end{equation}
On the other hand, since $P_{11}$ has full rank, then equality holds in (\ref{e:span_orthogonal_basis}). Thus, the first part of the claim is proved.

The second part can be shown in a similar fashion. $\Box$\\
\end{proof}

\section*{Acknowledgements}

The authors would like to thank the two anonymous reviewers for their comments and suggestions.

\bibliography{ofbm}

%\small

\small

\bigskip

\noindent \begin{tabular}{lcl}
Gustavo Didier & \hspace{3.5cm} & Vladas Pipiras \\
Mathematics Department& & Dept.\ of Statistics and Operations Research \\
Tulane University & & UNC-Chapel Hill \\
6823 St.\ Charles Avenue  & & CB\#3260, Smith Bldg. \\
New Orleans, LA 70118, USA & & Chapel Hill, NC 27599, USA \\
{\it gdidier@tulane.edu}& & {\it pipiras@email.unc.edu} \\
\end{tabular}\\

\smallskip

\end{document}